\documentclass[12pt,twoside]{article}
\usepackage{amsmath}
\usepackage{amsfonts}
\usepackage{amsthm}
\usepackage{amssymb}
\usepackage{graphicx}
\usepackage{multicol}

\begin{document}
\title{{\normalsize 
{\bf Non-trivialization probability of arc system in three-dimensional space}}} 
\author{{\footnotesize Akio Kawauchi}\\ 
\date{} 
{\footnotesize{\it Osaka Central Advanced Mathematical Institute, Osaka Metropolitan University}}\\  
{\footnotesize{\it Sugimoto, Sumiyoshi-ku, Osaka 558-8585, Japan}}\\  
{\footnotesize{\it kawauchi@omu.ac.jp}}} 
\date\,  
\maketitle 
\vspace{0.1in} 
\baselineskip=9pt 
\newtheorem{Theorem}{Theorem}[section] 
\newtheorem{Conjecture}[Theorem]{Conjecture} 
\newtheorem{Lemma}[Theorem]{Lemma} 
\newtheorem{Sublemma}[Theorem]{Sublemma} 
\newtheorem{Proposition}[Theorem]{Proposition} 
\newtheorem{Corollary}[Theorem]{Corollary} 
\newtheorem{Claim}[Theorem]{Claim} 
\newtheorem{Definition}[Theorem]{Definition} 
\newtheorem{Example}[Theorem]{Example} 

\begin{abstract} The type-specific knotting probability of an arc diagram is earlier  
defined by using  chord diagrams of ribbon surface-links in 4D space. 
By  modifying this notion, Non-Trivialization probability (simply  NT probability) 
for the arc diagram is introduced  and generalized to an arc system diagram. 
Some properties of the NT probability are shown. The  method  of transforming 
a polygonal arc in 3D space into a unique arc diagram  up to  isomorphisms  
earlier developed is generalized to a  polygonal arc system in 3D space to define 
the NT probability. 

\phantom{x}



\end{abstract}

\baselineskip=15pt

\bigskip

\noindent{\bf 1. Introduction}

\begin{figure}[hbtp]
\begin{center}
\includegraphics[width=14cm, height=8cm]{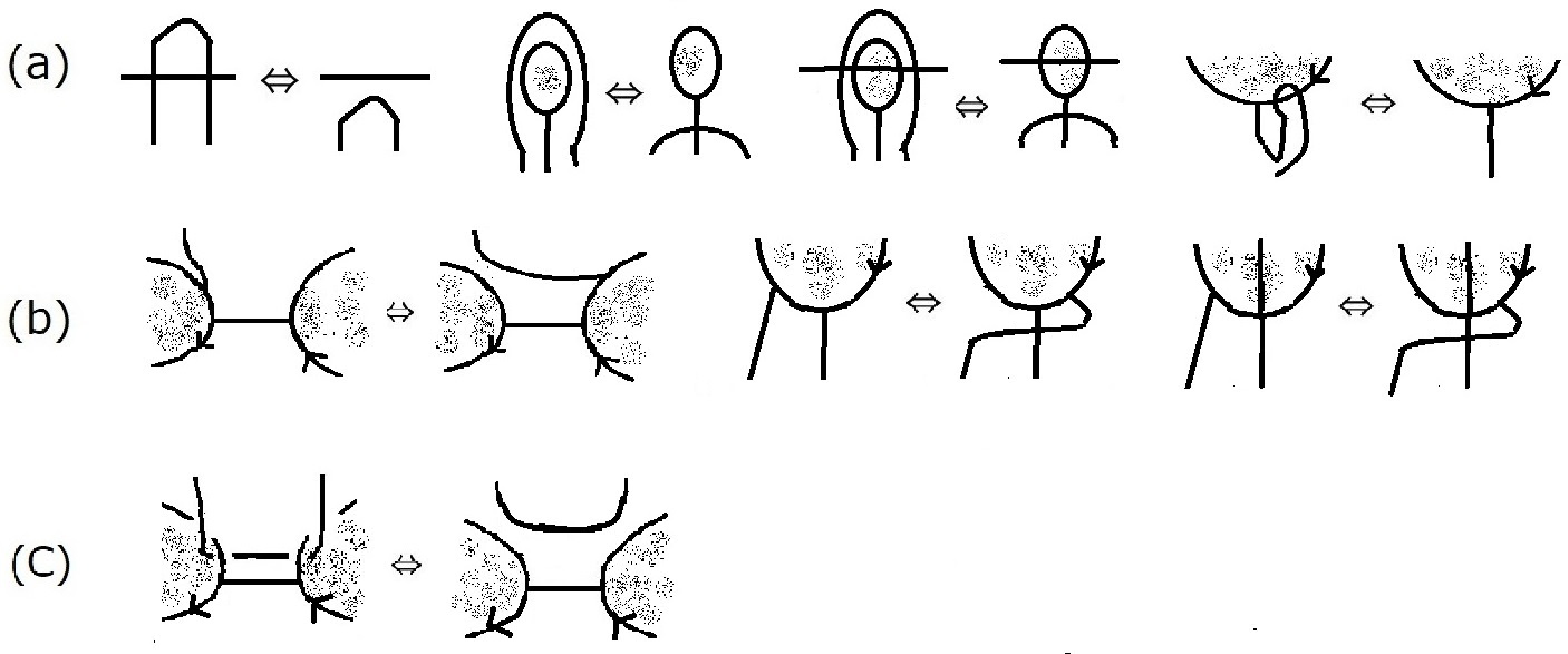}
\end{center}
\caption{A chord-homotopy move (a), a chord-loop move (b) and 
a fusion-fission chord move (c)}
\label{fig:bequiv}
\end{figure}

An $n$-{\it arc system} in 3-dimensional  space (simply 3D space) $R^3$ is 
the  system $L$ of mutually disjoint polygonal oriented arcs $a_i\, (i=1,2, \dots, n)$ in $R^3$ (Where unnecessary, the orientation of the arc $a_i$ is omitted). 
An arc system $D$ in an oriented plane $P$ is an $n$-{\it arc diagram}  if $D$ has only transversely meeting crossing points with upper-lower relation 
(apart from the end points)  in $P$.  
A standard arc system $D$ in an oriented plane $P$
is obtained as the image $\lambda(L)$ of an $n$-arc system $L$ in $R^3$ 
under an  orthogonal projection $\lambda: R^3\to P$ 
such that the singularity  ${\cal S}(L)$, the set of a point $x\in L$ with 
$|\lambda^{-1}(\lambda(x))|\geq 2$, consists of  a non-vertex point $x$ of 
an edge in $L$ with  $|\lambda^{-1}(\lambda(x))|=2$.
An {\it arc system diagram} is  an $n$-{\it arc diagram} for some $n$. 
A {\it chord graph} is a trivalent graph $(o; \alpha)$ in the 3D space $R^3$ consisting 
of a trivial oriented link $o$ (called a {\it based loop system}) and attaching 
simple arcs $\alpha$ (called a {\it chord system}), and 
a {\it chord diagram} is a planar representation $C(o;\alpha)$ of a chord graph 
$(o;\alpha)$ in a plane $P$, which represents a unique ribbon surface-link 
$F(o; \alpha)$  in 4-dimensional space (simply 4D space) $R^4$, \cite{[1], [2]}. 
Two chord diagrams $C=C(o,\alpha)$ and $C'=C(o',\alpha')$ in a plane $P$  are 
{\it isomorphic} if there is an orientation-preserving self-homeomorphism $f:P\to P$ 
sending $C$ to $C'$ which  preserves the orientations of the based loops $o$ and $o'$, 
and {\it equivalent} if the ribbon surface-link $F(o; \alpha)$ is equivalent to the ribbon 
surface-link $F(o'; \alpha')$. 
Equivalence between $C$ and $C'$ in $P$ is understood by  
a finite number of  the  moves $M_i\,(i=0,1,2)$  on chord diagrams, \cite{[3]}. 
This notion is also explained in Section~2.
The chord diagram $C=C(D)=C(o,\alpha)$ of an $n$-arc diagram $D$ is 
{\it basically equivalent} to the chord diagram $C'=C(D')=C(o',\alpha')$ of an 
$n$-arc diagram  $D'$ if,  up to isomorphisms of $C$, the based loop $o$  
is identical to the based loop $o'$ and the chord system $\alpha'$ is obtained from 
the chord system $\alpha$ by a finite number of  moves consisting of a 
{\it chord-homotopy move} (which is a plane homotopy move of a chord 
which does not change around the starting and terminal based loop system 
$o^s(a_i), o^t(a_i)\, (i=1,2,\dots,n)$), 
a {\it chord-loop move} and  a {\it fusion-fission chord move}, 
illustrated in Fig.~\ref{fig:bequiv}. 
These moves are cosequences of  the moves $M_i\,(i=0,1,2)$. 
An $n$-arc diagram $D$ is {\it basically equivalent} to  an $n$-arc diagram $D'$ if  
$C(D)$ is basically equivalent to  $C(D')$. 

\begin{figure}[hbtp]
\begin{center}
\includegraphics[width=8.5 cm, height=5 cm]{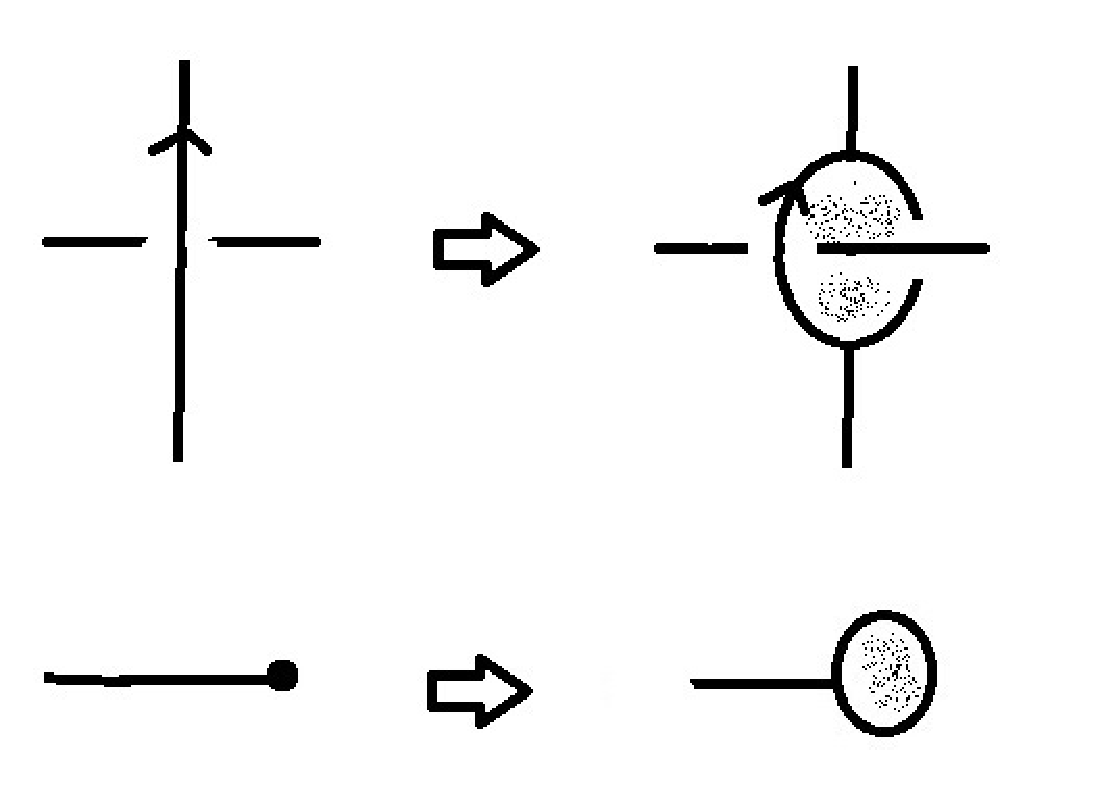}
\end{center}
\caption{Transforming the areas around a crossing point and an endpoint 
into chord diagrams}
\label{fig:trans}
\end{figure}

There is a method to transform every $n$-arc diagram $D$ of $c$ crossings  
into a chord diagram $C=C(D)$ with $c+2n$ based loops uniquely up to isomorphisms, 
where the trivial based loops  $o^s(a_i)$ and $o^t(a_i)$ are  attached to the 
starting point $s(a_i)$ and the terminal point $t(a_i)$ of every arc $a_i$ in $D$, 
respectively.   
This transformation was used  in \cite{[4]} and illustrated in Fig.~\ref{fig:trans}. 
The original $n$-arc diagram $D$ is recovered from the chord diagram 
$C=C(D)=C(o,\alpha)$ by taking the upper arc part of every based loop system $o$.
The first purpose of  this paper is to show the following theorem. 

\phantom{x}

\noindent{\bf Theorem~1.1.} For any $n$-arc diagram $D$ and any integer $d$ 
satisfying $1\leq d\leq 4$, there exists a basic equivalence invariant $\kappa_d(D)$ on 
$D$ that takes a rational value in the interval $[0,1]$ and is independent of choices of orientations of $a_i\, (i=1,2,\dots,n)$.

\phantom{x}

The proof of Theorem~1.1 is done in Section~3. This invariant $\kappa_d(D)$ is called 
 {\it Non-Trivialization probability} (simply {\it NT probability}) of {\it grade} $d$ 
 for the $n$-arc diagram $D$. 
The NT probabilities $\kappa_1(D)$ and  $\kappa_4(D)$ of grades $1$ and $4$ 
are also called the {\it  reduced NT probability} and the {\it   full NT probability} 
for $D$, respectively. 
The {\it  NT probability} $\kappa(D)$ for  an $n$-arc diagram $D$ is  the quadruplet  
\[\kappa(D)=(\kappa_1(D),\kappa_2(D), \kappa_3(D), \kappa_4(D)).\]
Some properties of the NT probability $\kappa(D)$ are shown in Theorem~3.1. 
The {\it mirror image} $D^*$ of an $n$-arc diagram $D$ is an $n$-arc diagram 
obtained from $D$ by changing the upper-lower relation of every crossing point of $D$. 
The NT probability $\kappa(D^*)$ for $D^*$ is generally different from $\kappa(D)$.  
The NT probability $\kappa(D)$
can be understood as the {\it NT probability of $D$ toward the future}, and 
$\kappa(D^*)$ as the {\it NT probability of $D$ toward the past}.
The second purpose of this paper is to show the following theorem.

\phantom{x}

\noindent{\bf Theorem~1.2.} For every $n$-arc system $L$  in $R^3$ up to 
orientation-preserving Euclidean isometries of $R^3$, there exists a unique system 
of $n$-arc diagrams $D(L)_{ij}\,(i,j=1,2,\dots,n)$ up to isomorphisms. 

\phantom{x}

The proof of Theorem~1.2 is done in Section~5.
By Theorem~1.1, the NT probability 
\[\kappa(D(L)_{ij})=(\kappa_1(D(L)_{ij}),\kappa_2(D(L)_{ij}),
\kappa_3(D(L)_{ij}),\kappa_4(D(L)_{ij}))\] 
of the $n$-arc diagrams $D(L)_{ij}$ is defined for every $i,j$. 
By taking the additive average 
\[\kappa_d(L)=\sum_{i,j=1}^n \kappa_d(D(L)_{ij})/n^2,\]
the {\it NT probability} 
$\kappa(L)=(\kappa_1(L),\kappa_2(L),\kappa_3(L),\kappa_4(L))$
of every $n$-arc system $L$ up to orientation-preserving Euclidean isometries 
of $R^3$ is defined. 

In Section~2, basic concepts on the chord diagram of an $n$-arc diagram are 
discussed. In Section~3, the NT probability of an $n$-arc diagram is defined 
and some properties of the NT probability are shown. 
In Section~4, the unique system of $n$-arc diagrams for an $n$-arc system in 
3D space $R^3$ is shown. In Section~5, computation examples of  $2$-arc systems 
in 3D space are given.

At the end of the introduction, the author refers to several historical notes. Arcs in 
$R^3$ from the random knotting viewpoint  independent of  the present viewpoint 
were studied , \cite{[5]}, \lq\lq KnotProt\rq\rq (https://knotprot.cent.uw.edu.pl/). 
This research come from the question: {\it How a linear scientific object such as 
a non-circular molecule (that is, a non-circular DNA, protein, linear polymer, etc.) 
is considered as a mathematical knot object?}
Also, a knotting probability on classical knots was  studied from the viewpoint of 
a random knotting, \cite{[6], [7]}. 
A knotting probability of an arc in $R^3$ was also  studied  by the author from 
a knotting structure of a  graph in $R^3$ but with the demerit that it depends 
on the heights of the crossing points of a  diagram of the  arc, 
\cite{[8], [9]}. The knotting probability of an arc diagram that resolves 
this demerit has been studied in \cite{[4]}. When attempting to generalize 
this knotting probability to the case of an $n$-arc diagram, we encountered 
a serious problem that does not manifest in the case of an arc, which motivated 
us to define the current NT probability.
The study of an arc system in 3D space is a natural generalization of the study 
of an arc in 3D space, \cite{[10]}. 

\phantom{x}

\begin{figure}[hbtp]
\begin{center}
\includegraphics[width=3cm, height=2cm]{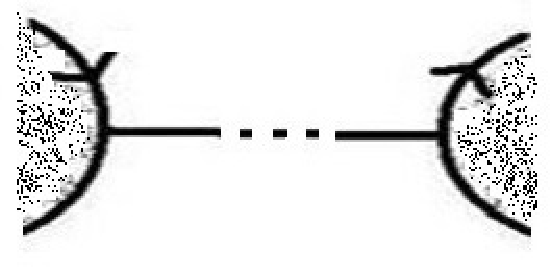}
\end{center}
\caption{A chord in an oriented regular chord diagram}
\label{fig:Orichord}
\end{figure}

\noindent{\bf 2. Basic concepts on chord diagram of arc system diagram}

A chord diagram $C=C(o;\alpha)$ is {\it regular} if 
the based loops $o$  bound mutually disjoint disks in the plane $P$ meeting only with 
the chords transversely. A regular chord diagram is {\it oriented} if 
every chord is attached to the  based loops as in Fig.~\ref{fig:Orichord}. 
An {\it orientable chord diagram} is a chord diagram $D$ 
which becomes an oriented regular chord diagram after applying a finite number 
of the move $M_0$ to $D$. {\it Unless otherwise mentioned,  a chord diagram 
means an orientable regular chord diagram.}
There are three moves $M_0$, $M_1$ and $M_2$ on  chord diagrams, \cite{[1]}.

\phantom{x}

\begin{figure}[hbtp]
\begin{center}
\includegraphics[width=10cm, height=6cm]{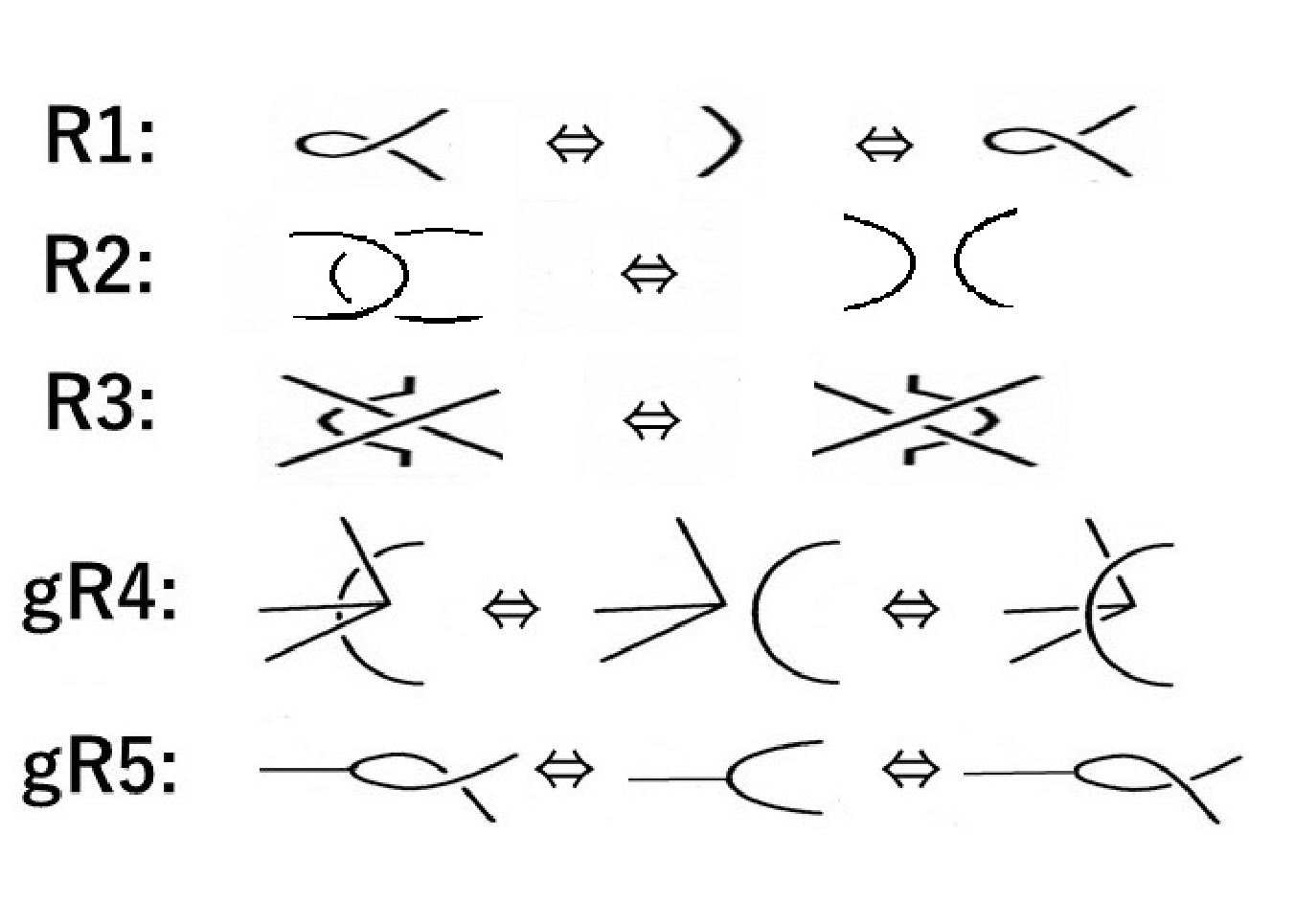}
\end{center}
\caption{Move $M_0$:\, Reidemeister moves $R_1$, 
$R_2$, $R_3$, $gR_4$, $gR_5$ for trivalent graph diagrams}
\label{fig:Rmoves}
\end{figure}

\noindent{\bf Move} ${\mathbf M_0.}$ This move consists of the Reidemeister 
moves $R_1$, $R_2$, $R_3$, $gR_4$, $gR_5$ as trivalent graphs in $R^3$, 
illustrated in Fig.~\ref{fig:Rmoves}, where note that  two arcs in the three arcs 
together with a vertex or an arc can be a part of a based loop although 
the orientation and the shadow of the based loop are omitted. 

\phantom{x}

\begin{figure}[hbtp]
\begin{center}
\includegraphics[width=8cm, height=2cm]{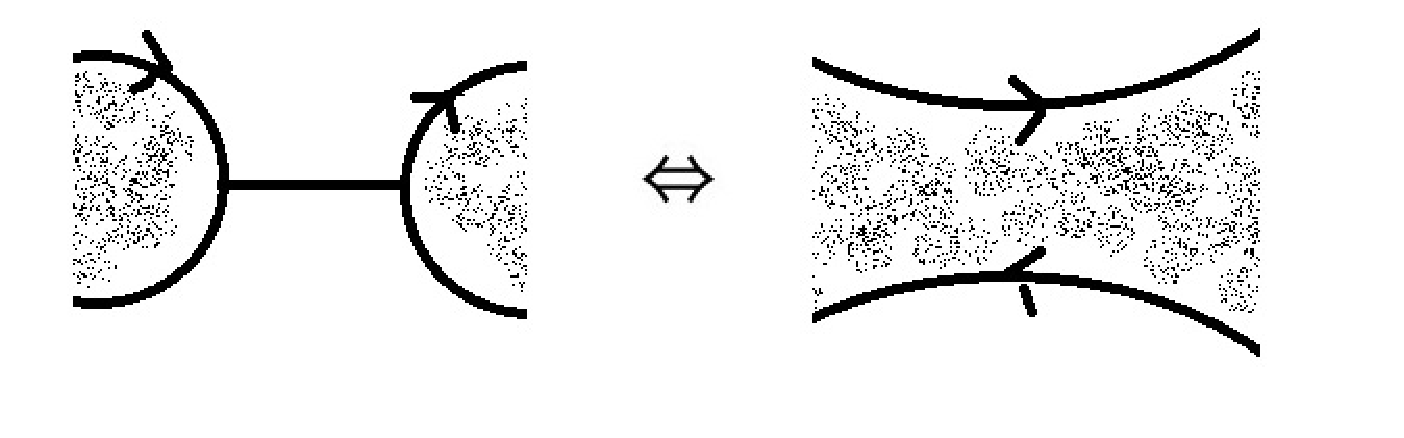}
\end{center}
\caption{Fusion-fission move $M_1$}
\label{fig:FF}
\end{figure}

\noindent{\bf Move} ${\mathbf M_1.}$ This move is the
{\it fusion-fission move}, illustrated in Fig.~\ref{fig:FF}, where the 
fusion operation is done only for a chord between different based loops. 

\phantom{x}

\begin{figure}[hbtp]
\begin{center}
\includegraphics[width=7cm, height=3cm]{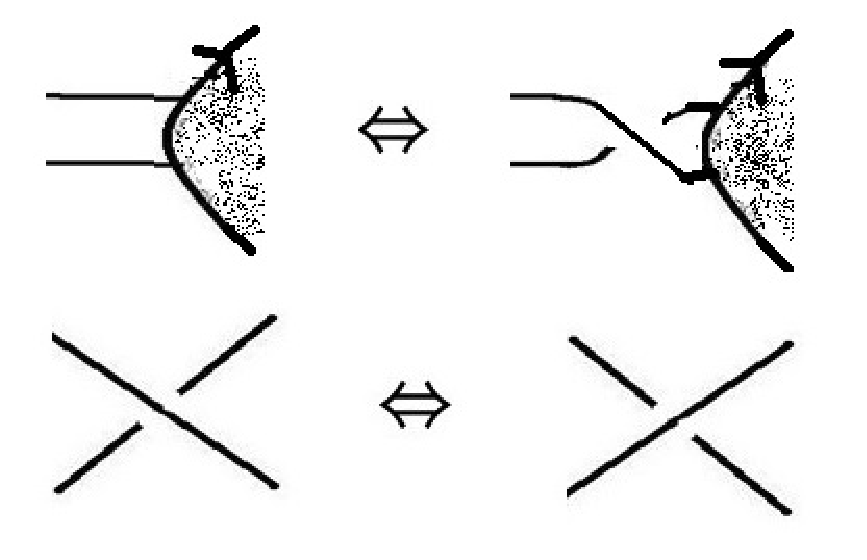}
\end{center}
\caption{Chord moves $M_2$}
\label{fig:chordmoves}
\end{figure}

\noindent{\bf Move} ${\mathbf M_2.}$  This move consists of moves on chords 
illustrated in Fig.~\ref{fig:chordmoves}. 

\phantom{x}

Two chord diagrams $C=C(o,\alpha)$ and $C'=C(o',\alpha')$ in a plane $P$  are 
{\it equivalent} if one is obtained from the other by a finite number of the moves 
$M_i\,(i=0,1,2)$ up to isomorphisms. 
The following result is known, \cite{[1], [3], [11]}. 

\phantom{x}

\noindent{\bf Lemma~2.1.}
Two ribbon surface-links $F(o;\alpha)$ and $F(o';\alpha')$ 
are equivalent (that is, sent by an 
orientation-preserving diffeomorphism of $R^4$)
if and only if the chord diagrams $C(o';\alpha')$ and  $C(o;\alpha)$ are equivalent.

\phantom{x}

\begin{figure}[hbtp]
\begin{center}
\includegraphics[width=8cm, height=6cm]{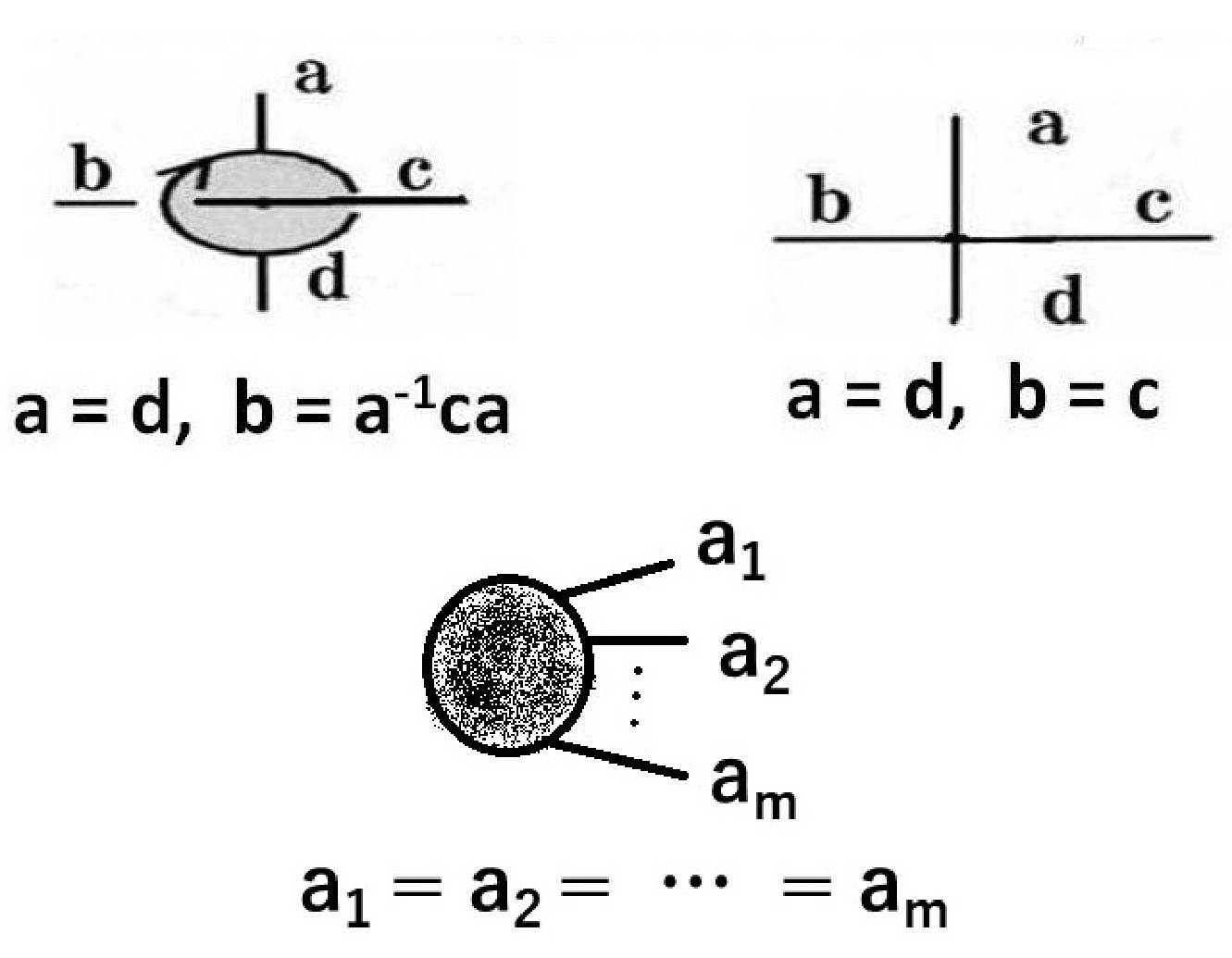}
\end{center}
\caption{Group relations}
\label{fig:grouprelations}
\end{figure}

The following lemma is basic to determine whether a chord diagram is trivial, 
that is, whether the surface-link is trivial,\cite{ [12]}. 

\phantom{x}

\noindent{\bf Lemma~2.2.} A  surface-link $F$ is a trivial surface-link 
in  $R^4$ if and only if the fundamental group $\pi_1(R^4\setminus F, x_0)$ is a free group with a meridian basis.

\phantom{x}

A useful invariant for determining whether a chord diagram $C=C(o;\alpha)$  
is trivial or not 
is the {\it finitely presented group} $\pi(C)$, which is 
identified with the fundamental group  $\pi_1(R^4\setminus F(C), x_0)$ 
and calculated as in Fig.~\ref{fig:grouprelations}, 
where the generator system is represented by 
a meridian system of $F(C)$ in $R^4$, \cite{[1], [13]}. 
Since all the edges of $C$ attaching to every based loop of $o$ have 
the same color in the group presentation $\pi(C)$, the based loop system $o$ 
is identified with the generator system of the finitely presented group $\pi(C)$.  
Note that the presented group  $\pi(C)$ is identical to the group of a virtual graph diagram interpreted from the chord diagram $C$, \cite{[1], [2]}.
For the group $\pi(C)$ of a chord diagram $C$, 
let  $\gamma: \pi(C)\to {\mathbf Z}$ be the epimorphism sending 
the generators of the presented group $\pi(C(D))$ to the unit $1\in {\mathbf Z}$, 
and $K(\gamma)$ the kernal subgroup  of $\gamma$.
The quotient group $K(\gamma)/K(\gamma)'$  for the commutator subgroup 
$K(\gamma)'$ of $K(\gamma)$  forms a finitely generated 
$\Lambda$-module $M(C)$, called the {\it module} of the chord diagram $C$, where 
$\Lambda$ denotes the integral Laurent polynomial ring ${\mathbf Z}[t, t^{-1}]$
identical to the integral group ring ${\mathbf Z}[{\mathbf Z}]$. 
The $\Lambda$-module $M(C)$ can be obtained from the group presentation of $\pi(C)$ 
by Fox's free differential calculus and is an invariant of an oriented 
chord diagram $C$ up to equivalences, \cite{[14], [15], [16]}. 
However, note that the $\Lambda$-module $M(C)$ depends on choices of orientations 
of the chord diagram $C$ in general (see Example~4.3 later).  

\phantom{x}

\noindent{\bf 3. Non-trivialization probability of arc system diagram}

Let $D$  be an $n$-arc  diagram of an arc system 
$L=\{a_i|\, i=1, 2,\dots, r\}$, and  $C(D)$  the chord  diagram of $D$. 
For an arc $a_i$ in $L$, let  
$c_i$ be the upper crossing number (i.e., the number of the upper crossing points) 
of $a_i$ in $D$.   
Let $C(a_i^D)$ be  the chord subdiagram  of the chord diagram $C(D)$ on the arc  
$a_i$ which has $c_i+2$ based loops including the two trivial based loops $o^s(a_i)$ and $o^t(a_i)$ of the endpoints of $a_i$. Note that the chord diagram $C(a_i^D)$ is obtained  from the chord diagram of the 1-arc diagram of $a_i$ in $D$ by adding 
some extra trivial based loops. 
The proof of Theorem~1.1 is done as follows.

\phantom{x}

\begin{figure}[hbtp]
\begin{center}
\includegraphics[width=11cm, height=7cm]{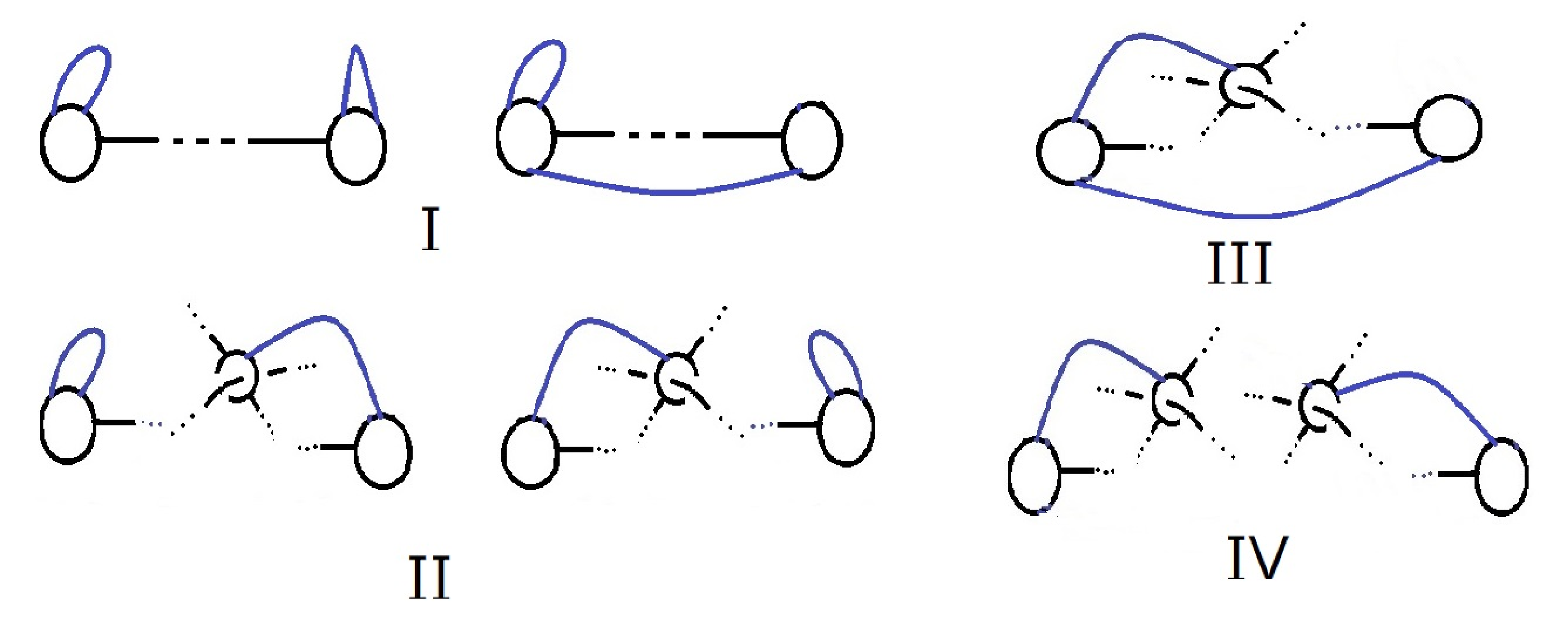}
\end{center}
\caption{Chord diagrams of types I, II, III, IV}
\label{fig:types}
\end{figure}

\noindent{\it 3.1: Proof of Theorem~1.1.}
An {\it $a_i$-based adjoint chord diagram} of the chord diagram $C(D)$ is a 
chord diagram  obtained from  $C(D)$ by connecting $o^s(a_i)$ and $o^t(a_i)$ 
to any different based loops of $C(a_i^D)$ with a single chord each not passing 
the other based loops of $C(D)$. 
Let $A(D, a_i)$ be the the system  of $a_i$-based adjoint chord diagrams of  
$C(D)$ consisting of $(c_i+2)^2$  chord diagrams.  
By removing some natural overlaps up to equivalences, the system $A(D, a_i)$ 
is reduced to the system  $A^R(D, a_i)$ consisting of $(c_i+1)^2+1$  
chord diagrams, which is classified into the  four  systems  
$A^{\tiny\mbox{I}}(D, a_i)$ of type I,  $A^{\tiny\mbox{II}}(D, a_i)$ 
of type II, $A^{\tiny\mbox{III}}(D, a_i)$ of type III, and  $A^{\tiny\mbox{IV}}(D, a_i)$ 
of  type IV, which are illustrated in Fig.~\ref{fig:types}, \cite{[4]}.
  
\phantom{x}

\noindent{\bf Type I system $A^{\tiny\mbox{I}}(D, a_i)$.} 
Here are two members of $a_i$-based adjoint chord diagrams 
$A_{ij}=A^{\tiny\mbox{I}}(D, a_i)_j\,(j=1, 2)$ of $C(D)$. One is 
the $a_i$-based adjoint chord diagram with two self-attaching additional chords. 
The other  is the $a_i$-based adjoint chord diagram with a self-attaching additional 
chord on $o^s(a_i)$ and an additional chord joining $o^s(a_i)$ with $o^t(a_i)$.  
Let $m_i^{\tiny\mbox{I}}$ be the number of 
non-trivial $a_i$-based adjoint chord diagrams in  $A^{\tiny\mbox{I}}(D, a_i)$. 

\phantom{x}

\noindent{\bf Type II system $A^{\tiny\mbox{II}}(D, a_i)$.} 
Here are  $2c_i$ members of  $a_i$-based adjoint chord diagrams 
$A_{ij}=A^{\tiny\mbox{II}}(D, a_i)_j\,(j=1,2,\dots, 2c_i)$ of $C(D)$, which are given by the additional chord pairs consist of a self-attaching additional chord on 
$o^s(a_i)$ (or $o^t(a_i)$, respectively) and an additional chord joining 
$o^t(a_i)$ (or $o^s(a_i)$, respectively) with a based loop except for $o^s(a_i)$ and 
$o^t(a_i)$. 
Let $m_i^{\tiny\mbox{II}}$ be the number of 
non-trivial $a_i$-based adjoint chord diagrams in  $A^{\tiny\mbox{II}}(D, a_i)$. 

\phantom{x}

\noindent{\bf Type III  system $A^{\tiny\mbox{III}}(D, a_i)$.} 
Here are  $c_i$ members of  
$a_i$-based adjoint chord diagrams $A_{ij}=A^{\tiny\mbox{III}}(D, a_i)_j\,(j =1, 2, \dots,n_i)$ 
of $C(D)$   where 
the additional chord pairs consist of an additional chord joining  $o^s(a_i)$ with  
$o^t(a_i)$  and an  additional chord joining $o^s(a_i)$ with  
a based loop except for  $o^s(a_i)$  and  $o^t(a_i)$. 
Let $m_i^{\tiny\mbox{III}}$ be the number of 
non-trivial $a_i$-based adjoint chord diagrams in  $A^{\tiny\mbox{III}}(D, a_i)$. 

\phantom{x}

\noindent{\bf Type IV system $A^{\tiny\mbox{IV}}(D, a_i)$.} 
Here are  $c_i(c_i - 1)$ members of 
$a_i$-based adjoint chord diagrams $A_{ij}=A^{\tiny\mbox{IV}}(D, a_i)_j\,
(j=1, 2, \dots, c_i(c_i -1))$ of $C(D)$  where 
the additional chord pair joins the pair of $o^s(a_i)$ and $o^t(a_i)$ 
with a pair of distinct based loops except for $o^s(a_i)$ and $o^t(a_i)$. 
Let $m_i^{\tiny\mbox{IV}}$ be the number of 
non-trivial $a_i$-based adjoint chord diagrams in  $A^{\tiny\mbox{IV}}(D, a_i)$. 

\phantom{x}

Since $2+2c_i+c_i+c_i(c_i-1)=1+(c_i+1)^2$,  the reduced  system  $A^R(D, a_i)$ 
consists of  just $1+(c_i+1)^2$ members of $a_i$-based adjoint chord diagrams. 
A general idea of NT probability is to measure to what extent non-trivial 
chord diagrams are included within the reduced $a_i$-based adjoint 
chord system  $A^R(D, a_i)$. 
By Lemma~2.2, an $a_i$-based adjoint chord diagram $A_{ij}$ in 
$A^R(D, a_i)$  is trivial if and only if the finitely presented group 
$\pi(A_{ij})$ of $A_{ij}$ is a free group with a basis in the generator system. 
For an $n$-arc diagram $D$ and an integer $d$ with $1\leq d\leq 4$, 
the {\it  $a_i$-based NT probability $\kappa_d(D, a_i)$ of grade $d$}  is defined as follows. 
\begin{eqnarray*}
\kappa_1(D, a_i)&=&m_i^{\tiny\mbox{I}}/2, \\
\kappa_2(D, a_i)&=&(m_i^{\tiny\mbox{I}}+m_i^{\tiny\mbox{II}})/(2+2c_i), \\
\kappa_3(D, a_i)&=&(m_i^{\tiny\mbox{I}}+m_i^{\tiny\mbox{II}}+m_i^{\tiny\mbox{III}})
/(2+3c_i),\\
\kappa_4(D, a_i)&=&
(m_i^{\tiny\mbox{I}}+m_i^{\tiny\mbox{II}}+m_i^{\tiny\mbox{III}}+m_i^{\tiny\mbox{IV}})
/(1+(c_i+1)^2).
\end{eqnarray*}
By Lemma 2.1 and definition, the value $\kappa(D)$ is  independent of choices 
of orientations of $a_i\,(i=1,2,\dots,n)$ and  invariant under  basic equivalences 
of an $n$-arc diagram $D$. This completes the proof of Theorem~1.1.

\phantom{x}

The {\it  NT probability of the $n$-arc diagram $D$ of grade $d$}  is defined by 
\[\kappa_d(D)=\sum_{i=1}^n\kappa_d(D, a_i)/n.\] 
Here,  $\kappa_1(D)$ and  $\kappa_4(D)$ are also called the 
{\it  reduced NT probability} and the {\it   full NT probability} of the $n$-arc diagram 
$D$, respectively. 
The {\it  NT probability} $\kappa(D)$ of  an $n$-arc diagram $D$ is defined by 
the quadruplet  
\[\kappa(D)=(\kappa_1(D),\kappa_2(D), \kappa_3(D), \kappa_4(D)).\]
The following notation is used. 

\phantom{x}

\noindent{\bf Notation.} The quadruplet 
$\kappa=(\kappa_1,\kappa_2,\kappa_3,\kappa_4)$ of  numbers 
$\kappa_i\,(i=1,2,3,4)$ in the interval $[0,1]$ is written as
$\kappa=1, \kappa> 0$, or $\kappa=0$ according to whether 
$\kappa_i =1\, (i=1,2,3,4)$, 
$\kappa_i >0\, (i=1,2,3,4)$, or $\kappa_i =0\, (i=1,2,3,4)$, respectively. 

\phantom{x}

An $n$-arc diagram $D$ of an $n$-arc system $L=\{a_i|\, i=1, 2,\dots, r\}$  is 
{\it inbound} if  $D$ is included in an $n$-component link diagram $\mbox{cl}(D)$ 
without extra crossing points and the endpoints $s(a_i), t(a_i)$ 
of $a_i\,(i=1,2,\dots,n)$ are in the same region in the rigions of the plane $P$ 
divided by $D$,  whose region is called  the  {\it front region}.  
The link diagram $\mbox{cl}(D)$ 
in $P$ is called the {\it closed link diagram} of  an inbound $n$-arc diagram $D$, 
which is the union of $D$ and a disjoint simple arc system  
$\gamma_i\,(i=1,2,\dots,n)$ with 
$\partial \gamma_i=\{s(a_i),  t(a_i)\}\,(i=1,2,\dots,n)$.  
The  following theorem on some properties of NT probability is proved.

\phantom{x}

\noindent{\bf Theorem~3.1.} The  NT probability $\kappa(D)$ of 
an $n$-arc  diagram $D$ of an arc system $L=\{a_i|\, i=1, 2,\dots, n\}$ 
has the following properties (1)-(3)

\phantom{x}

\noindent{(1)} $\kappa(D)=0$ if and only if   the chord diagram $C(a_i^D)$ is a trivial 
chord diagram for all $i$ and $C(D)$ is basically equivalent to a 
split sum of the trivial chord diagrams $C(a_i^D)$ for all $i$.

\phantom{x}

\noindent{(2)} $\kappa(D)=1$ if and only if  the quotient group  obtained 
from the presented group $\pi(C(D))$ by identifying each of
the generators $o^s(a_i)$ and $o^t(a_i)$ with a single generator 
in the based loop system of $C(a_i^D)$ 
is not a free group with a basis in the generator system  for every $i\, (i=1,2,\dots,n)$.

\phantom{x}

\noindent{(3)} If $D$ is an inbound $n$-arc  diagram and $D^*$ is the mirror image 
of $D$, then $\kappa(D)=\kappa(D^*)$ and $\kappa_1(D)= m/2n$ for an integer $m$ 
with $0 \leq m \leq 2n$ excluding $m=1, 3\, (n=2)$ and  $m=1\, (n \neq 2)$. 

\phantom{x}

A realization of the reduced NT probability $\kappa_1(D)$ of an inbound $n$-arc 
diagram  $D$ excluding the exceptional values in Theorem~3.1(3) 
is achieved in Example~4.4.

\phantom{x}

\noindent{\bf Proof of Theorem~3.1.} 
For (1), if $\kappa(D)=0$, then the chord diagram $C(D)$ is basically equivalent 
to a split sum of the chord diagrams $C(a_i^D)\, (i=1,2,\dots,n)$ and 
$\kappa(C(a_i^D))=0$ for every $i$. 
In particular, $\mbox{cl}_u(D)$ is a trivial link diagram since $\kappa(C(a_i^D))=0$ implies that the knot diagram  $\mbox{cl}_u(a_i^D)$ is a trivial knot diagram for every $i$, \cite{[4]}. The converse is obvious. This shows (1). 
For (2), let $A_{ij}$  be any $a_i$-based adjoint chord diagram of $C(D)$ in 
$A^{\Gamma}(D, a_i)$.  
The presented group $\pi(A_{ij})$ is obtained from the presented group $\pi(C(D))$ 
by identifying the generator  $o^t(a_i)$ with a based loop generator of $C(a_i^D)$ 
and the generator  $o^t(a_i)$ with a based loop generator of $C(a_i^D)$.  
By Lemma~2.2,  $\kappa(D)=1$ if and only if $\pi(A_{ij})$ is not a free group 
with based loop basis for all $i$ and $j$, completing the proof of (2). 
For (3), by Artin's spiinning construction of an inbound arc diagram $D$, 
the ribbon $S^2$-link $F(C)$ of the chord diagram $C=C(D)$ is sent 
to the ribbon $S^2$-link $F(C^*)$ of the mirror image $C^*=C(D^*)$ by 
an orientation-reversing diffeomorphism of $R^4$, \cite{[4]}.  This means that
$\kappa(D)=\kappa(D^*)$. To research the value of 
the reduced NT probability $\kappa_1(D)$ of an inbound arc diagram $D$, 
let $C=C(D)$ be the chord diagram of an inbound $n$-arc diagram $D$.
Let $\mbox{cl}(D)*v$ be a diagram of an $n$-leafed bouquet $L^+_v$ with vertex $v$ 
in the 3-sphere $R^3$  
obtained from the link diagram $\mbox{cl}(D)$ by adding $n$ simple arcs 
$\gamma'_i\, (i=1,2,\dots,n)$ disjoint except for the vertex $v$ such that 
$\gamma'_i$ joins a point of $\gamma_i$ to the vertex $v$ 
and $\mbox{cl}(D)*v$ has no crossing points other than the crossing points of 
$\mbox{cl}(D)$.
Let $\mbox{cl}(D)*v\setminus \mbox{Int}\gamma'_i$ is the subdiagram of 
$\mbox{cl}(D)*v$ which is a diagram of the subgraph  
$L^+_v\setminus \mbox{Int}\gamma'_i$ in $R^3$ 
obtained from   $L^+_v$  by removing the interior $\mbox{Int}\gamma'_i$ 
of the arc $\gamma'_i$. 
As explained in \cite[p.204]{[15]} and \cite{[17]}, by Artin's spinning construction, 
the fundamental group $G$ of 
$L^+_v\setminus \mbox{Int}\gamma'_i$ in $R^3$ is identical to the finitely presented group 
$\pi(C)$ of the chord diagram $C$ such that the generator set of $\pi(C)$ is 
a meridian system of  $L^+_v\setminus\{v_i\}$, and the fundamental grpup $G_i$ of 
$L^+_v\setminus \mbox{Int}\gamma'_i$ in $R^3$ is identical to the finitely presented group  
$\pi(C\cup\gamma_i)$ of the chord diagram $C(D)\cup\gamma_i$ such that 
the generator set of $\pi(C\cup\gamma_i)$ is a meridian system of  $G_i$.
If $\kappa_1(C)=1/2n$, then $G$ and $G_i$ for every $i$ except for one index $i$, 
say $i=1$ are free groups with basis in the specified generator systems.  
By applying Dehn's lemma 
to the compact exterior of $L^+_v\setminus \mbox{Int}\gamma'_i$ in $R^3$, 
the loop $a_i\cup\gamma_i$ bounds a disk $d_i$ in $R^3$ disjoint from 
the other component of $L^+_v\setminus \mbox{Int}\gamma'_i$.  
The disks $d_i\, (i=2,3,\dots, n)$ are made disjoint. 
Then the loop $a_i\cup\gamma_i$ also bounds a disk $d_n$ in $R^3$ disjoint from 
the other component of $L^+_v\setminus \mbox{Int}\gamma'_n$ and the disks 
$d_i\, (i=2,3,\dots, n)$. This means that the surface-links $F(C))$ and 
$F(C\cup\gamma_i)$ are trivial surface-links for all $i$, so that  $\kappa_1(C)=0$, 
a contradiction.  Hence $\kappa_1(C)\ne 1/2n$. 
Let $n=2$ and suppose that $\kappa_1(C)=3/4$.
Then the group $G_1$ or $G_2$, say $G_1$ must be  a free group with basis 
in the specified generator system. By  Dehn's lemma, the group $G_2$ is also 
a free group with basis in the specified generator system. Then $\kappa_1(C)=0$, 
a contradiction. Hence $\kappa_1(C)\ne 3/4$, showing (3). 
This completes the proof of Theorem~3.1. 

\phantom{x}

The following corollary gives a  criterion for an $n$-arc diagram $D$ with $\kappa(D)=1$

\phantom{x}

\noindent{\bf Corollary~3.2.} Let $M(C)$ be the module of the chord diagram $C=C(D)$ 
of an $n$-arc diagram $D$ of an arc system $L=\{a_i|\, i=1,2,\dots,n\}$.  Unless 
the module $M(C)$ 
is $\Lambda$-isomorphic to the direct sum $\Lambda^{n-1}\oplus TM(C)$ 
of a free $\Lambda$-module 
$\Lambda^{n-1}$ of rank $n-1$ and a   $\Lambda$-tortion module 
$TM(C)$ generated by two $\Lambda$-elements, the NT probability $\kappa(D)=1$. 

\phantom{x}

\noindent{\bf Proof of Corollary~3.2.} By Theorem~3.1 (2), $\kappa(D)=1$ if 
the $\Lambda$-module $M(A_{ij})$  is not  isomorphic to the free $\Lambda$-module 
$\Lambda^{n-1}$ of rank $n-1$ for all $i, j$. 
Note that the module $M(C)$ has  $\Lambda$-rank $n-1$ and 
the $\Lambda$-module $M(A_{ij})$  is a quotient 
$\Lambda$-module of $M(C)$ by two $\Lambda$-elements. 
If the $\Lambda$-module $M(A_{ij})$ for some $i, j$ is a free module of rank $n-1$, 
then the module $M(C)$ must be  $\Lambda$-isomorphic to 
the direct sum $\Lambda^{n-1}\oplus TM(C)$ for a  
$\Lambda$-tortion module $TM(C)$ generated by two $\Lambda$-elements,
which contradicts the assumption of the module $M(C)$. Thus, $\kappa(D)=1$. 
This completes the proof of Corollary~3.2.

\phantom{x}

For a 1-arc diagram $D$, the  {\it knotting probabilities of types} I, II, III, 
and IV denoted by $p^{\tiny\mbox{I}}(D)$, $p^{\tiny\mbox{II}}(D)$,  
$p^{\tiny\mbox{III}}(D)$, and $p^{\tiny\mbox{IV}}(D)$, 
respectively were defined in \cite{[4]}:
\[
p^{\tiny\mbox{I}}(D)=m_1^{\tiny\mbox{I}}/2, \quad
p^{\tiny\mbox{II}}(D)=m_1^{\tiny\mbox{II}}/2c_1, \quad
p^{\tiny\mbox{III}}(D)=m_1^{\tiny\mbox{III}}/c_1, \quad
p^{\tiny\mbox{IV}}(D)=m_1^{\tiny\mbox{IV}}/c_1(c_1-1).
\]
The {\it type-specific knotting probability} $p(D)$ of a 1-arc diagram $D$ 
is defined to be the quadruplet 
\[p(D)=(p^{\tiny\mbox{I}}(D), p^{\tiny\mbox{II}}(D), p^{\tiny\mbox{III}}(D), 
p^{\tiny\mbox{IV}}(D)).\]
The inconvenience of this definition is that the  
knotting probabilities $p^{\tiny\mbox{II}}(D)$, $p^{\tiny\mbox{III}}(D)$,
$ p^{\tiny\mbox{IV}}(D)$ 
are not defined when $c_1=0$, and the knotting probability $p^{\tiny\mbox{IV}}(D)$ 
is not defined when $c_1=1$, although these values were defined by $0$, respectively.
This inconvenience is serious to consider a general  $n$-arc diagram $D$ because 
$c_i=0$ happens often. From this reason,  the NT probability $\kappa(D)$ 
of an $n$-arc diagram $D$ is introduced. 

\phantom{x}

\noindent{\bf 4. Computation example on $n$-arc diagram }

\begin{figure}[hbt]
\begin{center}
\includegraphics[width=6 cm,height=2 cm]{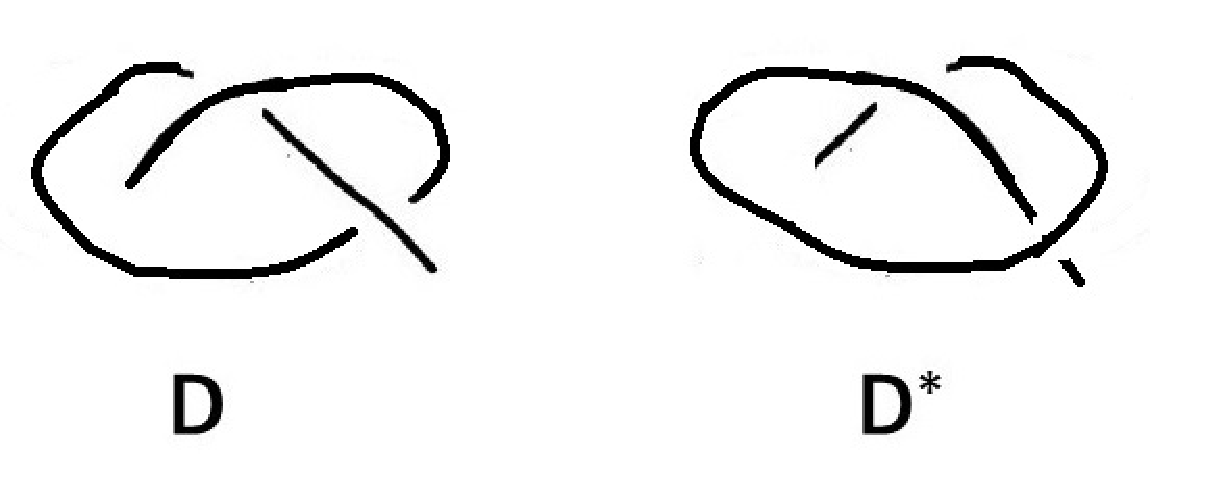}
\end{center}
\caption{A 1-arc diagram $D$ with 2 crossings and its mirror image $D^*$}
\label{fig:1A2CD}
\end{figure}

\begin{figure}[hbt]
\begin{center}
\includegraphics[width=14 cm,height=7cm]{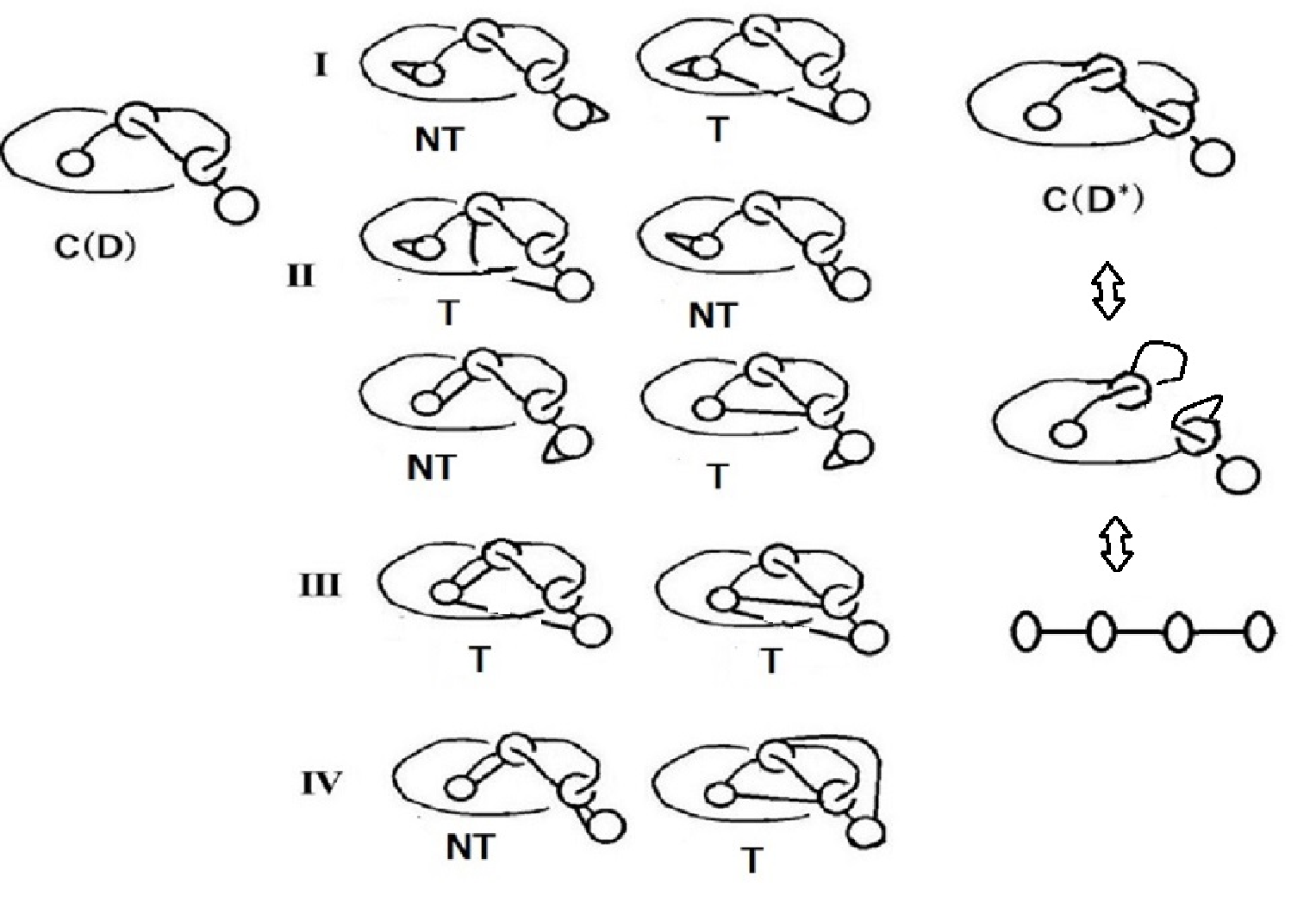}
\end{center}
\caption{The system $A^R(D,a_1)$ of  $a_1$-based adjoint chord diagrams of 
the chord diagram $C(D)$ of the diagram $D$ in  Fig.~\ref{fig:1A2CD}  
and a basic equivalence of the  mirror image $D^*$ of $D$ to a chord diagram 
without crossing}
\label{fig:1A2CDCal}
\end{figure}

This section presents several examples of how to calculate the NT probability for 
$n$-arc diagrams. In the examples, the notation T or NT attached to 
each chord diagram in the system $A^R(D, a_i)$ means  a trivial chord diagram 
or a nontrivial chord diagram, respectively. 
If the upper crossing number $c_1$ of  a 1-arc diagram $D$ is a previous known 
integer greater than or equal to $2$, then  the type-specific knotting probability 
$p(D)$ and the  NT probability $\kappa(D)$ are computable from each other.  
Some examples of the knotting probability $p(D)$ for the 1-arc diagram $D$  are 
calculated, \cite{[4]}. So, the NT probabilities $\kappa(D)$ of thse examples are obtained 
by rewriting the values of  these examples of $p(D)$. 
In the following example, the knotting probability was calculated in \cite{[4]}, 
but here the NT probability is directly calculated.

\phantom{x}

\noindent{\bf Example~4.1.}
For the 1-arc diagram $D$ with $c_1=2$ in 
Fig.~\ref{fig:1A2CD},  
the system $A^R(D, a_1)$ of $a_1$-based adjoint chord diagrams on the chord system 
$C=C(D)$ of $D$ is listed in Fig.~\ref{fig:1A2CDCal}.
The knotting probability $p(D)$ of $D$ is caluculated  as 
$p^{\tiny\mbox{I}}(D)=1/2$,  $p^{\tiny\mbox{II}}(D)=2/4=1/2$, 
$p^{\tiny\mbox{III}}(D)=0/2=0$, and $p^{\tiny\mbox{IV}}(D)=1/2$. 
Hence $p(D)=(1/2,1/2,0,1/2)$. The NT probability $\kappa(D)=\kappa(D,a_1)$ of $D$ 
is calculated as $\kappa_1(D)=1/2$, 
$\kappa_2(D)=(1+2)/(2+4)=1/2$, $\kappa_3(D)=(1+2+0)/(2+4+2)=3/8$, and 
$\kappa_4(D)=(1+2+0+1)/(2+4+2+2)=2/5$. Hence 
$\kappa(D)=(1/2, 1/2, 3/8, 2/5)$. 
The mirror image  $D^*$ of $D$ in Fig.~\ref{fig:1A2CD} is basically equivalent 
to a chord without crossing, so that $\kappa(D^*)=0$ by Theorem~1.1.
This shows that the NT probability 
$\kappa(D^*)$ of the mirror image  $D^*$ of an $n$-arc diagram $D$ is 
generally different from the NT probability $\kappa(D)$ of $D$.  

\phantom{x}

\begin{figure}[hbtp]
\begin{center}
\includegraphics[width=11cm, height=3cm]{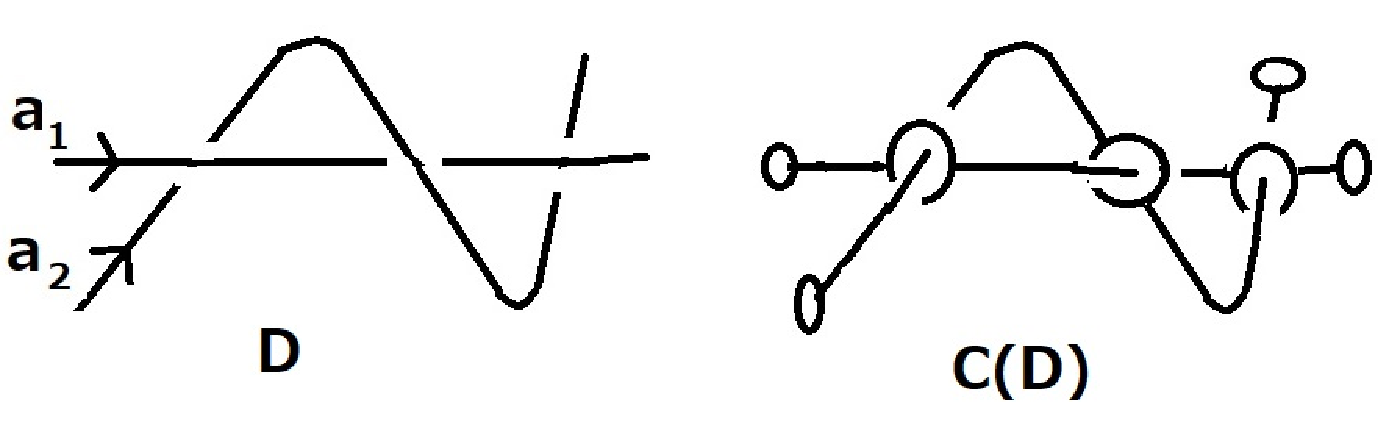}
\end{center}
\caption{A 2-arc diagram $D$ with 3 crossings and the chord diagram $C(D)$}
\label{fig:2A3CD}
\end{figure} 

\begin{figure}[hbtp]
\begin{center}
\includegraphics[width=16cm, height=8cm]{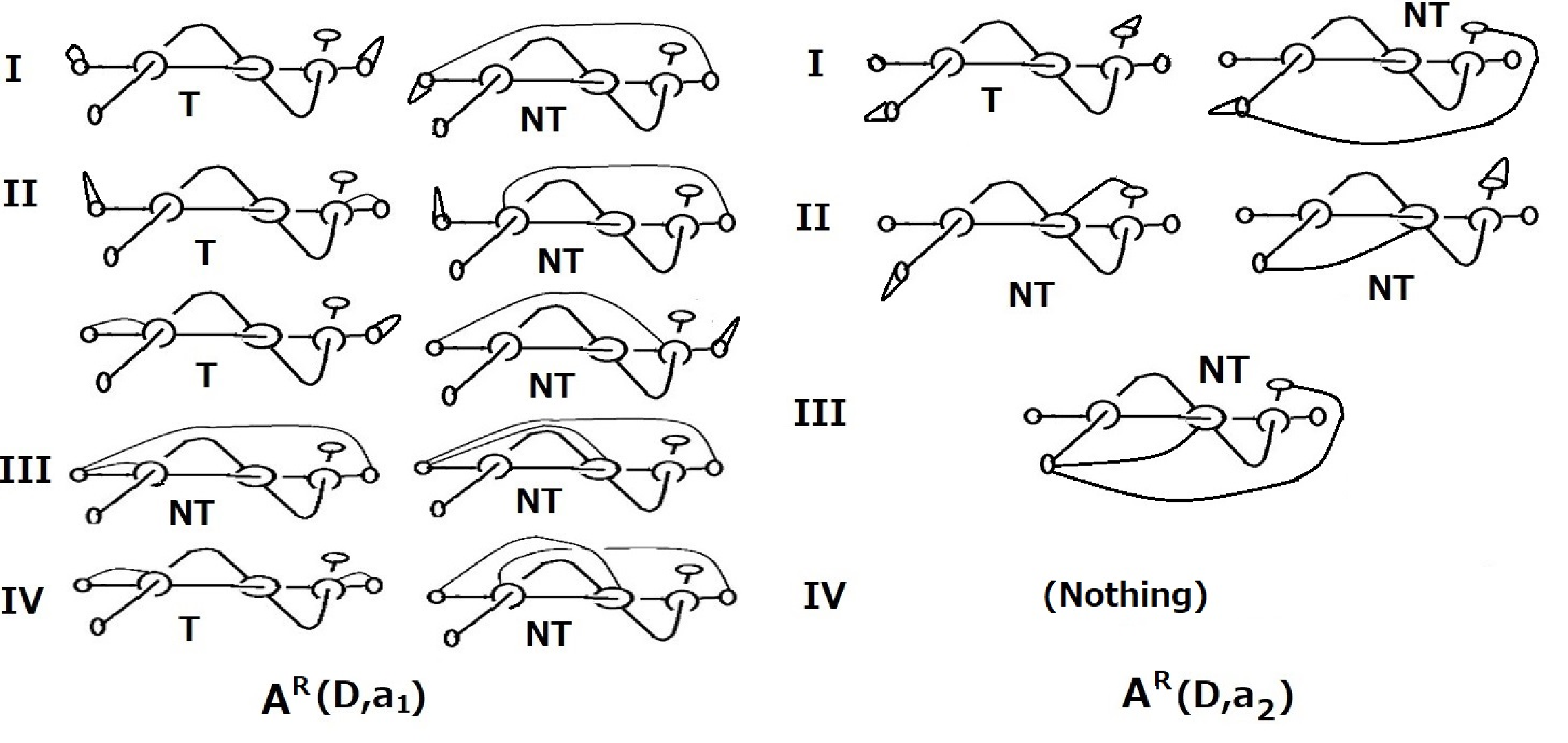}
\end{center}
\caption{The system $A^R(D,a_i)\, (i=1,2)$ of  $a_i$-based adjoint chord diagrams of 
the chord diagram $C(D)$ of the diagram $D$ in  Fig.~\ref{fig:2A3CD}}
\label{fig:2A3CDCal}
\end{figure}

\noindent{\bf Example~4.2.} The 2-arc diagram $D$ with  
$c_1=2$ and $c_2=1$ in  Fig.~\ref{fig:2A3CD}  induces the 
system $A^R(D, a_i)$ of  $a_i$-based adjoint chord diagrams $(i=1,2)$  in  
Fig.~\ref{fig:2A3CDCal}. Then  
$\kappa_1(D,a_1)=1/2$, 
$\kappa_2(D,a_1)=(1+2)/(2+4)=1/2$, 
$\kappa_3(D,a_1)=(1+2+2)/(2+4+2)=5/8$, 
$\kappa_4(D,a_1)=(1+2+2+1)/(2+4+2+2)=3/5$. 
Hence, $\kappa(D,a_1)=(1/2, 1/2, 5/8, 3/5)$. 
Also, $\kappa_1(D,a_2)=1/2$, 
$\kappa_2(D,a_2)=(1+2)/(2+2)=3/4$, 
$\kappa_3(D,a_2)=(1+2+1)/(2+2+1)=4/5=\kappa_4(D,a_2)$. 
Hence, $\kappa(D,a_2)=(1/2, 3/4, 4/5, 4/5)$. 
Thus, $\kappa(D)=(1/2, 5/8, 57/80, 7/10)$. 
Although $D$ is not inbound, the mirror image $D^*$ of the 2-arc diagram $D$  is 
basically equivalent to $D$, so that by Theorem~1.1, 
\[\kappa(D^*)=\kappa(D)=(1/2, 5/8, 57/80, 7/10).\]

\phantom{x}

The following example is an example on Theorem~3.1 (2) and Corollary~3.2.

\phantom{x}

\begin{figure}[hbt]
\centerline{\includegraphics[width=10 cm, height=4cm]{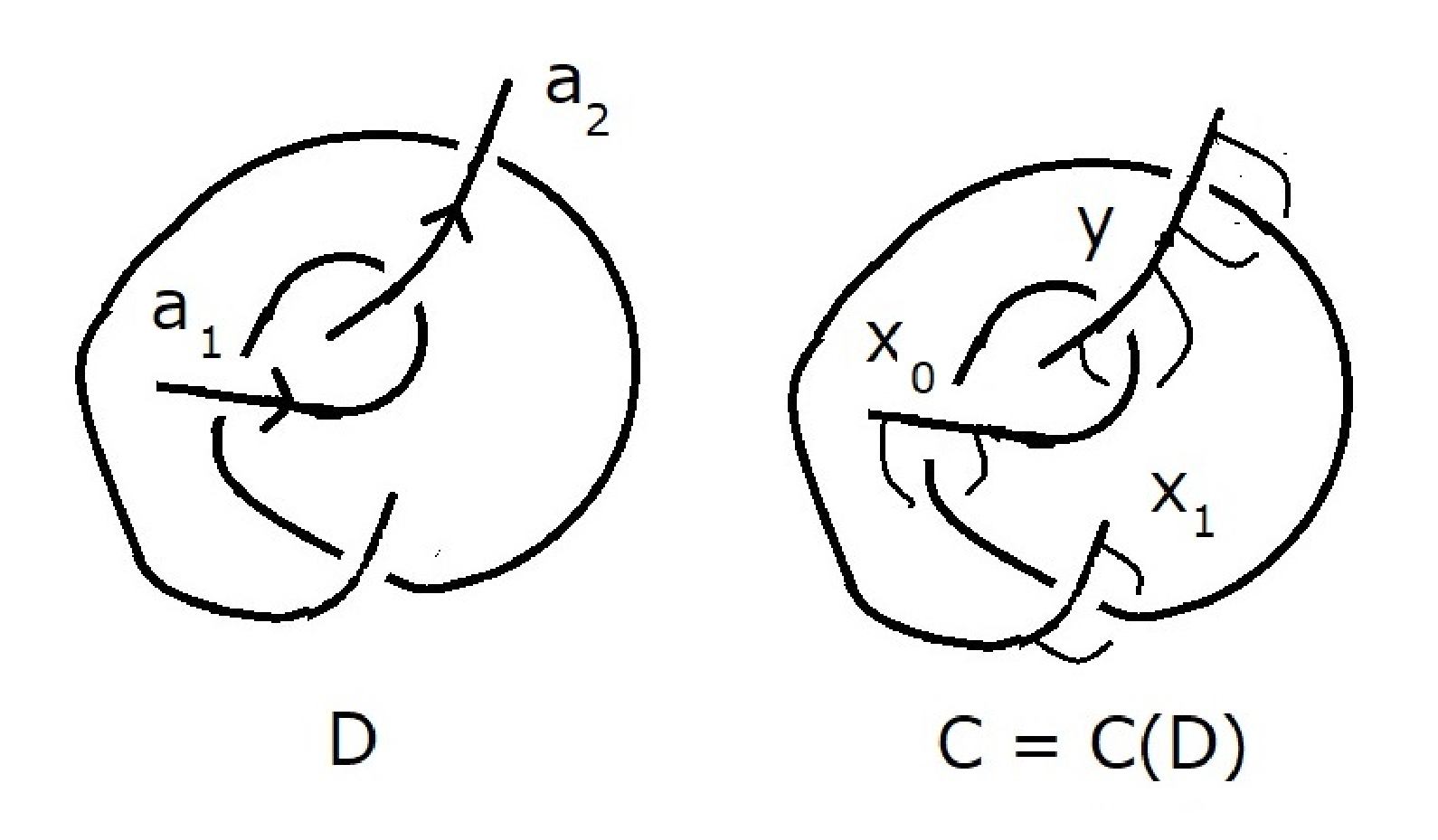}}
\caption{A $2$-arc  diagram $D$ and the chord diagram $C=C(D)$}
\label{fig:NFDiagram}
\end{figure}

\noindent{\bf Example~4.3.}
The finite presented group $\pi(C)$ of the chord diagram $C=C(D)$ of 
the 2-arc diagram $D$, illustratexd in Fig.~\ref{fig:NFDiagram} is given by 
\[\pi(C) = <x_0, x_1, y|\, x_0^{-1}x_1^{-1}yx_1y^{-1}x_1x_0y^{-1}x_0^{-1}y>.\]
The module $M(C)$ is calcutated to be $\Lambda$-isomorphic to the 
$\Lambda$-ideal $(t-2, 2t-1)$ of $\Lambda$, which is a torsion-free 
$\Lambda$-module of rank one, but not a free $\Lambda$-module. Thus, 
$\kappa(D)=1$ by Corollary~3.2.  Let $D'$ be the 2-arc diagram obtained from $D$ 
by changing the orientation of  only the arc $a_2$ into the opposite orientation, 
and $C'=C(D')$ the chord diagram of $D'$. 
Then the group $\pi(C')$ has the presentation
\[\pi(C') = <x_0, x_1, y|\, x_0^{-1}x_1^{-1}y^{-1} x_1yx_1x_0yx_0^{-1}y^{-1}>.\]
The module $M(D')$ is  calcutated to be $\Lambda$-isomorphic to the direct sum 
$\Lambda\oplus \Lambda/(t^2-t+1)$, which  
cannot be used to determine that $\kappa(D')=\kappa(D)=1$.

\phantom{x}

The following example shows that the value $\kappa_1(D)$ 
of an inbound $n$-arc diagram $D$ excluding the exceptional values 
in Theorem~3.1(3)  is complete.

\phantom{x}

\begin{figure}[hbt]
\centerline{\includegraphics[width=13.5 cm,height=4cm]{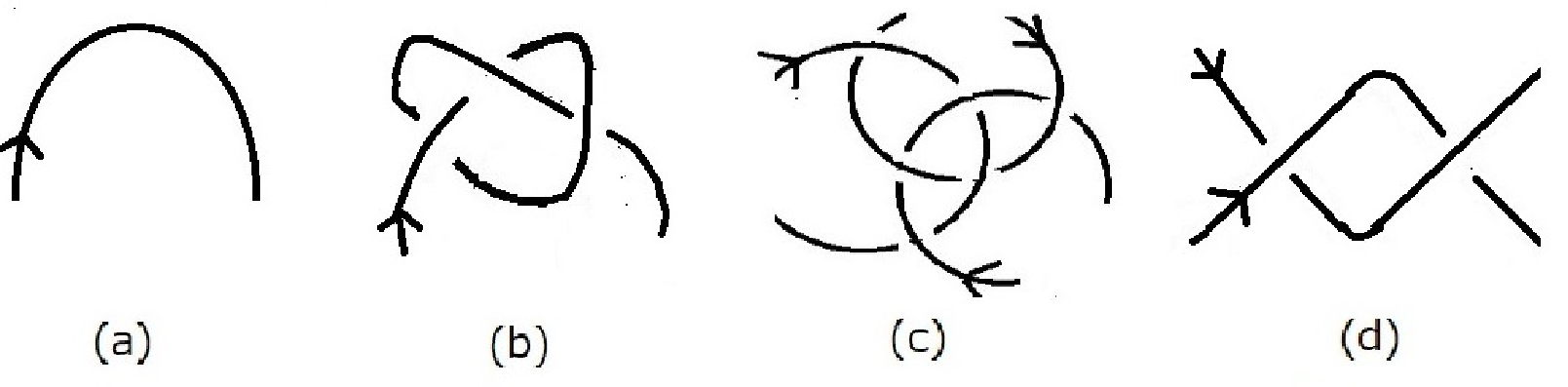}}
\caption{ Diagrams (a), (b), (c), (d) }
\label{fig:Reduced}
\end{figure}

\noindent{\bf Example~4.4.} By Theorem~3.1 (3), the reduced NT probability 
$\kappa_1(D)$ of an inbound $n$-arc diagram $D$ of $L=\{a_i|\, i=1,2,\dots,n\}$ 
takes a value $m/2n$ with $0\leq m\leq 2n$ and $m\ne 1, 2n-1$.  
Any such value $m/2n$ is  realized by an  inbound $n$-arc  diagram $D$. 
The split union diagram $D$ of $n$ copies of diagram $D_1$ of 
Fig.~\ref{fig:Reduced} (a) has $\kappa_1(D)=0$ by definition.
The split union diagram $D$ of $n$ copies of the diagram $D_2$ of  
Fig.~\ref{fig:Reduced} (b) has $\kappa_1(D)=1$ since $\kappa_1(D_2)=1$, \cite{[4]}.  
So, assume that $2\leq m\leq 2n-2$. If $m=2p$ for an integer $p$, 
then $1\leq p\leq n-1$. Then take as $D$ a split diagram of  $p$ copies 
of  the diagram $D_2$ of Fig.~\ref{fig:Reduced} (b) 
and $n-p$ copies of the diagram $D_1$. Then   $\kappa_1(D)=p/n=m/2n$.   
If $m=2p+3$ for an integer $p$, then  $0\leq p\leq n-3$ and $n\geq 3$. 
If $n=3$, then $p=0$ and take as $D$ the diagram $D_3$ of 
Fig.~\ref{fig:Reduced} (c)  where $\kappa_1(D)=1/2=m/2n$, 
computed as $\kappa_1(D,a_i)=1/2\, (i=1,2,3)$. 
If $n\geq 4$, then take as $D$ the split diagram 
of $p$ copies of the diagram $D_2$, one copy of $D_3$,  
and $n-3-p$ copies of the diagram $D_1$. 
Then $\kappa_1(D)=(1\times p+ (1/2)\times 3)/n=m/2n$. 
If a non-split diagram is wanted, then $D$ can be transformed 
into a non-split inbound diagam by using R2-moves in Reidemeister move 
(cf. Fig.~\ref{fig:Rmoves}), where the resulting diagram is basically equivalent 
to the original diagram and thus, the reduced NT probability is unchanged. 
For $n\geq 3$, it is shown that there is an inbound non-split $n$-arc diagram $D$ with 
$\kappa_1(D)=(2n-1)/2n$.  
Let $o_1\cup o_2$ be a split union of two trivial loops $o_1, o'_1$ in $R^3$. 
Let $\gamma'_1$ be a simple arc in $R^3$ joining a point of $o_1$ and a point of 
$o'_1$ such that the union $o_1\cup o'_1\cup \gamma'_1$ is 
a non-trivial (= non-planar) graph in $R^3$. 
Construct a trivial $(n-1)$-leafed bouquet 
$(\cup_{i=2}^n o_i)\cup(\cup_{i=2}^n \gamma'_i)$ 
with vertex $v$ in a solid torus neighborhood $N(o'_1)$ of $o'_1$ in $R^3$ 
such that $o_2=o'_1$ and the leaf $o_i (i\geq 3)$ is parallel to $o'_1$ in $N(o'_1)$, 
where  $\gamma'_i$ denotes a leg (i.e.,  a simple arc joining a point of $o_i$ to 
the vertex $v$ with 
$\gamma'_i\setminus\{v\}\,(i=2,3,\dots,n)$ disjoint). 
By sliding the endpoint of the arc $\gamma'_1$ in $o'_1$ to the vertex $v$ along 
$o_2\cup\gamma'_2$, regard  $\Gamma=(\cup_{i=1}^n o_i)\cup(\cup_{i=1}^n \gamma'_i)$ 
as an $n$-leafed bouquet in $R^3$.  Let $B'$ be a regular neighborhood 3-ball of 
$\cup_{i=1}^n \gamma'_i$ in $R^3$, and $D$ an inbound  $n$-arc diagram in a plane $P$ 
of a proper $n$-arc system $L=\mbox{cl}(\Gamma\setminus\Gamma\cap B') 
=\cup_{i=1}^n a_i$ in a 3-ball $B\subset R^3$ with 
$a_i=\mbox{cl}(o_i\setminus o_i\cap B')$.   
The $n$-leafed bouquet $\Gamma$ is a non-trivial graph in $R^3$, 
by Artin's spinning construction, the chord diagram $C=C(D)$ of the inbound  
$n$-arc diagram$D$ is non-trivial. 
Since $\Gamma$ is a non-trivial graph and 
$(\cup_{i=1}^n o_i)\cup \gamma'_1\cup\gamma'_j$ 
is a non-splittable non-trivial graph for $j\geq 2$, 
and $(\cup_{i=1}^n o_i)\cup(\cup_{i=2}^n \gamma'_i)$
is a trivial graph in $R^3$,  
it is computed via Artin's spinning construction that $\kappa_1(D,a_1)=1/2$ and 
$\kappa_1(D,a_i)=1$ for all $i\geq 2$, so that 
$\kappa_1(D)=(1/2+1\times (n-1))/n=(2n-1)/2n$, as desired. 

\phantom{x}

\noindent{\bf 5. Unique $n$-arc diagram system for $n$-arc system in 3D space}

Let $L=\{a_i|\, i=1,2,\dots,n\}$ be an oriented polygonal $n$-arc system  
in  3D space $R^3$.
In this section, first it will be observed how the projection image of  
$L$ to an oriented plane $P$ is approximated by an $n$-arc diagram. 
Next, it will be shown how a unique system of finitely many $n$-arc chord diagrams 
up to isomorphisms is constructed from any given $n$-arc system $L$ in $R^3$. 
Let
\[S^2=\{ u\in  R^3 |\, ||u||=1\}\]
be the unit sphere, where $||\phantom{u}||$ denotes the norm on $R^3$. 
Every element $u\in S^2$ is called a {\it unit vector from the origin $(0,0,0)$}.
For a unit vector $u\in R^3$, let $P(u)$ be the oriented plane containing 
the origin $(0,0,0)$ such that the unit vector $u$ is a positive normal vector to $P(u)$.  
For  a unit vector  $u=(a,b,c)$, the plane $P(u)$ denotes the plane $ax+by+cz=0$ 
in $R^3$ with  $xyz$-coordinate system.
The orthogonal projection from $R^3$ to the plane $P(u)$ is called the 
{\it projection along} the unit vector $u\in S^2$ and denoted by 
\[\lambda(u):R^3\to P(u).\] 
For a small positive number $\delta$, a $\delta$-{\it approximation} of 
the projection $\lambda(u):R^3\to   P(u)$ along $u\in S^2$  is the projection 
$\lambda(u'):R^3\to    P(u')$  
along a unit vector $u'\in S^2$  with $||u'-u||<\delta$, 
which is denoted by 
\[\lambda(u)^{\delta}: R^3\to P(u).\]
The following theorem is obtained.

\phantom{x}

\noindent{\bf Theorem~5.1.}  
Let $L$ be an oriented polygonal $n$-arc system in the 3-space $R^3$ with 
$xyz$ coordinate system, and  $\lambda=\lambda(u_z):R^3\to P$  the projection along the 
$z$-axis unit vector $u_z=(0,0,1)\in S^2$ to the plane $P=P(u_z)$. 
For any sufficiently small positive number $\delta$,  
the projection  $\lambda$ has a $\delta$-approximation
\[\lambda^{\delta}: R^3\to P\] 
such that the projection image $\lambda^{\delta}(L)=\lambda(u_z)^{\delta}(L)$ 
is an $n$-arc diagram  uniquely determined up to isomorphisms from  
$L$ and the projection $\lambda$. 

\phantom{x}

This proof is done by refining an argument transforming a classical knot 
in $R^3$ into a regular knot diagram and was done for a 1-arc in $R^3$, \cite{[10], [14]}.   
Let $L=\{a_i|\, i=1,2,\dots,n\}$ be an $n$-arc system in $R^3$, and  
$s(a_i)$ and $t(a_i)$ the starting point  and the terminal point of the arcs
$a_i\, (i=1,2,\dots, n)$, respectively. 
An {\it edge line} of $L$ is  the oriented straight line in $R^3$ extending an edge in $L$ 
with orientation induced from the oriented arc $a_i$ of $L$ containing the edge. 
A {\it front line} of $L$ is the oriented straight line  in $R^3$ containing the points 
$s(a_i)$ and $t(a_i)$ with  orientation from $s(a_i)$ to $t(a_i)$ for an arc  $a_i$ of $L$. 
For any plane $P$ in $R^3$, the {\it great circle} $C$  of $P$ in $R^3$ 
is  the great circle obtained by the intersection of $S^2$ and 
the plane $P(u)$ (containing the origin $0$) which is parallel to $P$.  
The {\it trace set} $T$ of  an $n$-arc system $L$ is the subset of $S^2$ consisting  
of  the great circles and the unit vectors obtained from  $L$ arising in the following 
cases  (i) and (ii):

\phantom{x}

\noindent{(i)} The great circle $C$ of  $S^2$  of the plane $P$ 
in $R^3$ determined from an edge  or  front line  $\ell$ of $L$ 
and a vertex $v$ of $L$ disjoint from $\ell$. 

\medskip

\noindent{(ii)} The unit vectors $\pm u_{\ell}$ of an edge line or a front line $\ell$ of $L$ 
and the unit vectors $\pm u_{\eta}\in S^2$  of  
a straight line $\eta$ that meets three lines $\ell_i\, (i=1,2,3)$ at three points, 
where $\ell_i\, (i=1,2,3)$ consist of edge or front lines such that no two lines lie 
on the same plane.

\phantom{x}

In (i), note that the trace set $T$ contains the following great circles: 
the great circle of the plane determined by two parallel distinct front or edge lines, and
the great circle of the plane determined by two distinct edge or front lines meeting 
at a point. In particular, the trace set $T$ contains the great circle of 
any plane containing 
any three distinct front or edge lines. This means that
the unit vectors $\pm u(\eta)\in S^2$  
of  a straight line $\eta$ in $R^3$ that meets three lines $\ell_i\, (i=1,2,3)$ 
at three points, where $\ell_i\, (i=1,2,3)$ consist of edge or front lines such that some 
two lines lie on the same plane. 
In (ii), note that a line $\eta$ meeting $\ell_i\, (i=1,2,3)$ at three points 
is unique by uniqueness of the solution to a system of three linear equations 
with regular matrix coefficients. Also, note that  if a unit vector $u\in S^2$ 
is in the trace set $T$, then the unit vector $-u$ is also in $T$. 
The following lemma is used for the proof of Theorem~5.1. 

\phantom{x}

\noindent{\bf Lemma~5.2.} For every unit vector $u\in S^2\setminus T$, 
the projection image $\lambda(u)(L)$ is an $n$-arc diagram in the plane $P(u)$.  
Further, the $n$-arc diagram  $\lambda(u)(L)$ up to isomorphisms is independent 
of  any choice of a unit vector $u'$ in the connected region $R(u)$ of $S^2\setminus T$ 
containing $u$. 

\phantom{x}

\noindent{\bf Proof of Lemma~5.2.}
If a unit vector $u\in S^2$ is not in (i), then every edge of $L$ and  every front line 
$\ell_{\gamma}$ are embedded into the plane $P(u)$ by the projection $\lambda(u)$. 
If a unit vector $u\in S^2$ is  in neither (i) nor (ii), then the set of vertices 
of $L$ is embedded into $P(u)$ by the projection $\lambda(u)$   whose image is disjoint 
from the image of any open edge of $L$. 
In particular, any two distinct parallel edge lines are disjointedly embedded into $P(u)$. 
Further, the images of the edges of $L$ meet only in the images of the open edges 
of $L$. Thus, if a unit vector $u\in S^2$ is in neither (i) nor (ii), namely if 
$u\in S^2\setminus T$, then the meeting points among the edges of $L$ 
consisting of double points between two open edges of $L$ and hence 
the projection image $\lambda(u)(L)$ is an $n$-arc diagram in the plane $P(u)$. 
The $n$-arc diagram  $\lambda(u')(L)$ is unchanged up to isomorphisms 
for any unit vector $u'$ in a connected open neighborhooh of $u$ in $S^2\setminus T$, 
so that the $n$-arc diagram  $\lambda(u')(L)$ is unchanged up to isomorphisms 
for any unit vector $u'$ in the connected region $R(u)$. 
This completes the proof of Lemma~5.2

\phantom{x}

For a number $x$ with $0\leq x\leq 1$, let $x^c=\sqrt{1-x^2}$. Then the proof of Theorem~5.1 is done as follows:

\phantom{x}

\noindent{\it 5.3: Proof of Theorem~5.1.}
The idea of the proof is to specify a unique connected region 
$R(u)$ of $S^2\setminus T$ adjacent to the unit vector $u_z=(0,0,1)\in S^2$. 
By taking a positive number $r$ sufficiently small, 
take a unit vector $u(r)=(r,0,r^c)\in S^2$   
so that the unit vector $u(x)=(x,0,x^c)$ with $0<x\leq r$ does not meet $T$ 
except for the great circle $x^2+z^2=1$ in $T$ (if it is in $T$). 
Then the connected region 
$R(u)$ is taken to be 
the unique connected region of $S^2\setminus T$ which is adjacent to the unit vector 
$u(r)$ and contains the unit vector 
$u(r, \theta)=(r\cos \theta, r\sin \theta, r^c)$ of $S^2$ 
for a sufficiently small positive number $\theta$, which is  orthogonal to the plane 
$P(r,\theta):\, (r\cos \theta)x+(r\sin \theta)y+r^c z=0$.  
The connected region $R(u)$ is uniquely specified. 
By Lemma~5.2, the image $\lambda(u)(L)$ of $L$ under the orthogonal projection 
$\lambda(u):R^3 \to P(u)$ is an $n$-arc diagram in the plane $P(u)$ 
which is a unique 
$n$-arc diagram up to isomorphisms independent of choices of a unit vector 
$u'$ in the connected region $R(u)$.  
Since the connected region $R(u)$ is adjacent to the unit vector $u_z\in S^2$,  
for every $\delta>0$  there is a unit vectors $u'\in R(u)$   
with $||u'-u||<\delta$, so that the orthogonal projection 
$\lambda(u'):R^3 \to P(u')$ is a desired $\delta$-approximation   
$\lambda^{\delta}: R^3\to P$. 
This completes the proof of Theorem~5.1.  

\phantom{x}  

Note that the vectors 
$u_X=(r^c\cos\theta,r^c\sin\theta,-r)$, 
$u_Y=(-\sin\theta,\cos\theta,0)$  and $u_Z=u(r,\theta)$ are 
mutually orthogonal unit vectors such that the plane $P(r,\theta)$ is spanned by 
$u_X, u_Y$. The following matrix $O(r,\theta)$ and the transpose 
matrix $O(r,\theta)^t=(u_X^t\, u_Y^t \, u_Z^t)$ (which are mutually 
inverse matrices of determinant one) are important for calculations. 
\[O(r,\theta)=\left(
\begin{array}{c}
u_X\\
u_Y\\
u_Z
\end{array}
\right)
=
\left(
\begin{array}{ccc}
r^c\cos\theta & r^c\sin\theta & -r\\
-\sin\theta & \cos\theta & 0\\
r\cos \theta & r\sin \theta & r^c
\end{array}
\right), \]

\[O(r,\theta)^t=(u_X^t\, u_Y^t\, u_Z^t) =
\left(
\begin{array}{ccc}
r^c\cos\theta &  -\sin\theta & r\cos\theta   \\
r^c\sin\theta & \cos\theta & r\sin \theta  \\
-r & 0 & r^c
\end{array}
\right).
\]
By the mutually orthogonal unit vectors $u_X$, $u_Y$  and $u_Z$, 
the new $XYZ$-coordinate system is considered. 
Since $u_X O(r,\theta)^t=(1\,0\,0)$, $u_Y O(r,\theta)^t=(0\,1\,0)$ and 
$u_Z O(r,\theta)^t=(0\,0\,1)$, the point $p=(a,b,c)$ in the $xyz$-coordinate system 
is transformed into the point $p=(a',b',c')$ in the $XYZ$-coordinate system 
by the matrix equation
\[(a'\, b'\, c')=(a\, b\, c) O(r,\theta)^t.\]
From now, the proof of Theorem~1.2 is done as follows.

\phantom{x}

\noindent{\it 5.4: Proof of Theorem~1.2.}
It is shown how a unique system of finitely many $n$-arc chord diagrams 
up to isomorphisms is constructed from any given $n$-arc system $L$ in $R^3$. 
For  every arc $a_i$ of $L$, let $\ell_i$ be the front line of $a_i$.   
For every arc $a_j$ of $L$ not belonging to $\ell_i$, let  
$e_j$ be the edge of $a_j$ which pops for the first time from  
$\ell_i$ when a point is going on $a_j$  along  the orientation of $a_j$, 
and $\ell_{ij}$ be the edge line of $e_j$.
Let  $P_{ij}$  be the plane in $R^3$ containing the origin $(0,0,0)$  which is 
parallel to the plane $P(\ell_i, \ell_{ij})$ spanned by the straight lines $\ell_i$ 
and $\ell_{ij}$ in $R^3$.  Then a new unique $XYZ$-coordinate system of $R^3$ is 
defined as follows. 
Take the $X$-axis given by the unit vector of $\ell_i$. 
The $Y$-axis is defined  by the normal vector of the oriented line $\ell^*_{ij}$ in 
$P(\ell_i, \ell_{ij})$ meeting $\ell_i$ orthogonally at $s(a_i)$ such  that the inner product of 
the unit vectors of $\ell^*_{ij}$ and $\ell_{ij}$ is positive. 
The $Z$-axis is constructed in $R^3$ by taking a positive unit vector $u_{ij}$ 
orthogonal to the oriented plane $P_{ij}$ with the origin $(0,0,0)$. 
By applying Theorem~5.1 to the image $\lambda(L)$ of $L$ under the projection 
$\lambda=\lambda(u_{ij}):R^3\to P_{ij}$, a unique $n$-arc diagram $D(L)_{ij}$ of $L$ 
up to isomorphisms is obtained. 
For every arc $a_j$ of $L$ belonging to $\ell_i$, let $P_{ij}$ be a plane  in $R^3$
containing the origin $(0,0,0)$  and the straight line $\ell_i$, and 
$D(L)_{ij}$ the union of disjoint  intervals $I_{ij}\, (j=1,2,\dots,n)$  in $P_{ij}$  
such that the interval $I_{ij}$ is obtained from  the image interval  
${\hat\lambda}_{ij}(a_j)$ of the arc $a_j$ 
in $\ell_i$ under the orthogonal projection ${\hat\lambda}_{ij}:R^3\to \ell_i$ 
by  pushing it  into $P_{ij}\setminus\ell_i$. 
This completes the proof of Theorem~1.2.

\phantom{x}

Note that for every arc $a_j$ of $L$ belonging to $\ell_i$, the $n$-arc diagram 
$D(L)_{ij}$ has no crossing point, so that the NT probability $\kappa(D(L)_{ij})=0$. 
An $n$-arc system $L$ in $R^3$ is {\it even} if
the front pop planes $P_{ij}\, (j=1,2,\dots, n_i)$ for every arc $a_i$ of $L$ 
are taken as the same plane $P$ in $R^3$ (although the orientations of 
$P_{ij}\, (j=1,2,\dots, n_i)$  do not necessarily coincide due to the definition of  $Y$-axis.) 
Then the plane $P$ is called the {\it front pop plane} of $L$. 
By the definition of an even $n$-arc system, note that 
for every $n$-arc system $L$ and a plane $P$ in $R^3$ 
such that the intersection $P\cap L$ is a finite point set containing 
the boundary point set $\partial L$ of $L$,  there is an even $n$-arc system 
$L^+$ in $R^3$ 
obtained from $L$ by adding some small edges in $P$ to $\partial L$  
such that all of the front pop planes of $L^+$ coincide with the  plane $P$. 
The following corollary is direct from the definition of NT probability of $L$. 

\phantom{x}

\noindent{\bf Corollary~5.5.} 
Assume that  the orthogonal projection image of an even $n$-arc system $L$ in 
$R^3$  into the front pop plane $P$ is an $n$-arc diagram  $D(L)$. 
Then each of the $n$-arc diagrams $D(L)_{ij}\,(i,j=1,2,\dots,n)$ 
of $L$ is isomorphic to $D(L)$ or the mirror image $D(L)^*$ of $D(L)$. 
If the $n$-arc diagram $D(L)$ of an even $n$-arc system $L$ in $R^3$ 
is inbound, then the NT probability $\kappa(L)$ of $L$ is equal to the NT probability $\kappa(D(L))$.

\phantom{x}

\noindent{\bf Proof of Corollary~5.5.} The first assertion is obtained from Theorem~1.2 
and definition.  The second assertion is obtained from Theorem~3.1 (3), completing the proof of Corollary~5.5.

\phantom{x}

In particular, for every  $n$-arc system $L$ embedded in a plane $P$ of $R^3$, 
the NT probability $\kappa(L)$ of $L$ is $0$.

\phantom{x}

\noindent{\bf 6. Computation example of  arc system in 3D space}

In this section, two elementary examples on a  2-arc system $L$  in $R^3$ 
are computed to show how to perform the calculation, where a point 
$(x,y,z)$ of $R^3$ is represented as $(x,y)^z$  in the plane $P$ defined by $z=0$  
and  $z$ is called the {\it height} of the point $(x,y)$. 

\phantom{x}

\begin{figure}[hbt]
\centerline{\includegraphics[width=6cm,height=4 cm]{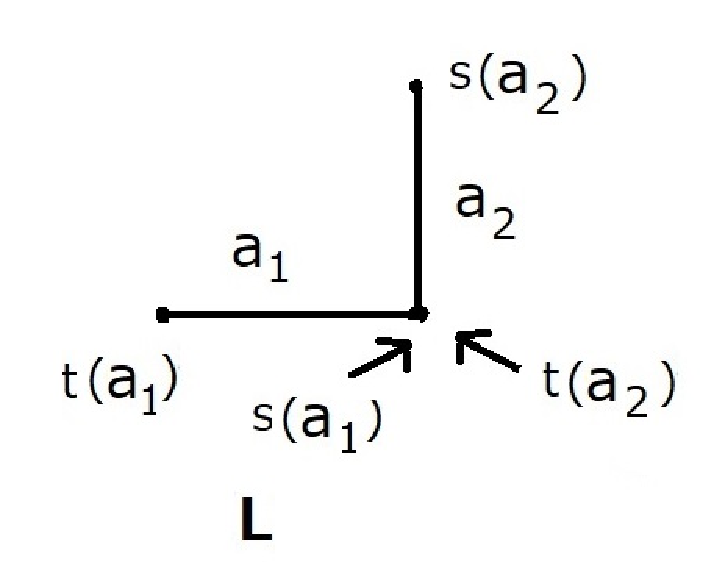}}
\caption{The non-even $2$-arc system $L=\{a_1, a_2\}$ represented in the 
$xyz$-coordinate system of $R^3$ with the following ordered vertices $(x,y)^z$:\,  
$a_1$: $s(a_1)=(0,0)^0$, $t(a_1)=(-1,0)^0$,\quad 
$a_2$: $s(a_2)=(0,1)^0$, $t(a_2)=(0,0)^1$}
\label{fig:L0arcxy_z}
\end{figure} 

\phantom{x}

\begin{figure}[hbt]
\centerline{\includegraphics[width=7cm,height=2.5 cm]{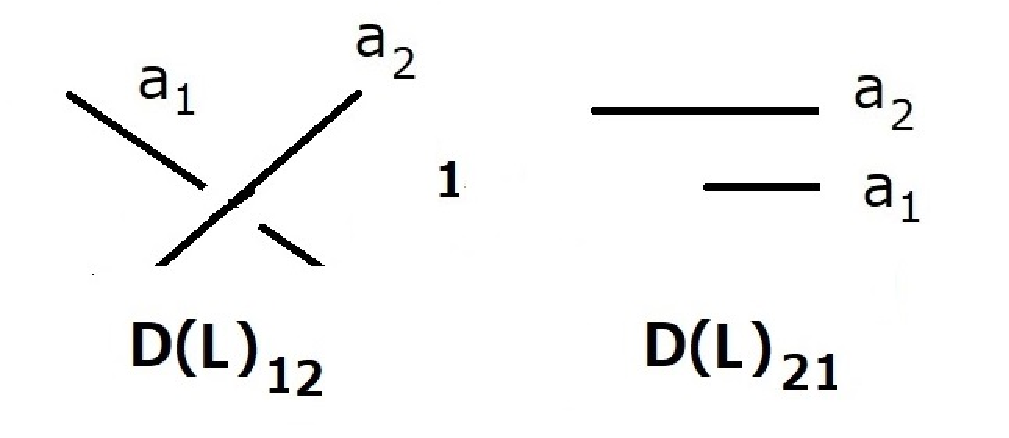}}
\caption{The 2-arc diagrams $D(L)_{12}$ and $D(L)_{21}$ of the 2-arc system 
$L$ in Fig.~\ref{fig:L0arcxy_z}}
\label{fig:D0arcxy_z}
\end{figure}

\noindent{\bf Example~6.1.} Let $L=\{a_1, a_2\}$ be a non-even 2-arc system in $R^3$ 
illustrated in  Fig.~\ref{fig:L0arcxy_z}. 
By definition,  $\kappa(D(L)_{11})=\kappa(D(L)_{22})=0$. To compute $\kappa(D(L)_{12})$, 
let $P_{12}$ be the front pop plane induced from $\partial a_1$ and $s(a_2)$, 
which coincides with the plane $P$ with standard orientation. 
The 2-arc system $L$ in the $XYZ$-coordinate system by the matrix 
$O(r,\theta)^t=(u_X^t\, u_Y^t \, u_Z^t)$ is given as follows. 
\[\begin{array}{ll}
a_1:& s(a_1) = (0,0)^0, \quad t(a_1)=(-r^c\cos\theta, \sin\theta)^{-r\cos\theta}, \\
a_2:& s(a_2)=(r^c\sin\theta, \cos\theta)^{ r\sin\theta}, \quad t(a_2)=(-r, 0)^{r^c}. 
\end{array}\]
By taking $r$ and $\theta$ to be sufficiently small positive numbers, 
the numbers $r'=r/r^c$ and  $\theta'=\tan\theta$ 
are taken as sufficiently small positive numbers. 
To obtain a 2-arc diagram $D(L)_{12}$ up to isomorphisms, 
one can use the $XYZ$-coordinate system obtained 
by applying the following matrix 
\[V(r,\theta)=((1/r^c\cos\theta)u_1^t\, (1/\cos\theta)u_2^t\,(1/r^c)u_3^t)=
\left(
\begin{array}{ccc}
1 & -\theta' & r' \\
\theta' & 1 & r'\theta'\\
-r'       & 0 & 1
\end{array}
\right)\]
instead of $O(r,\theta)^t=(u_1^t\, u_2^t \, u_3^t)$. 
Then the 2-arc system $L$ is given as follows. 
\[\begin{array}{ll}
a_1&: s(a_1) = (0,0)^0, \, t(a_1)=(-1, \theta')^{-r'},\\
a_2&: s(a_2)=(\theta', 1)^{r'\theta'},\, t(a_2)=(-r', 0)^1.  
\end{array}\]
Then the 2-arc diagram $D(L)_{12}$ is illustrated  in Fig.~\ref{fig:D0arcxy_z} and 
the computation results  
 $\kappa(D(L)_{12},a_1)=(1/2,1/2,1/2,1/2)$ and  $\kappa(D(L)_{12},a_2)=0$ are obtained.
Hence, $\kappa(D(L)_{12})=(1/4,1/4,1/4,1/4)$.   
Next the front pop plane $P_{21}$ is induced from $\partial a_2$ and $s(a_1)$, 
which is spanned by the unit vectors 
$u_X=(0,-1/\sqrt{2},1/\sqrt{2}), u_Y=(0,-1/\sqrt{2},-1/\sqrt{2})$
for the new $XY$-coordinate system of the plane $P_{21}$.  The $Z$-axis is defined by 
the unit vector $u_Z=(1,0,0)$. 
The plane $P_{21}$ is defined by $x=0$. 
Let $U$ be the $(2,2)$-matrix given by 
\[U=(u_X^t\, u_Y^t\,u_Z^t)=
\left(
\begin{array}{ccc}
0                &  0              & 1 \\
-1/\sqrt{2 } & -1/\sqrt{2} & 0 \\
1/\sqrt{2}    & -1/\sqrt{2} & 0
\end{array}
\right).\]
Since $u_X U=(1\,0\,0)$, $u_Y U=(0\,1\,0)$ and $u_Z U=(0\,0\,1)$, 
the 2-arc system $L$ in the $XYZ$-coordinate system 
of $R^3$ is obtained by applying the matrix $U$. However, 
since it is sufficient to obtain a 2-arc diagram up to isomorphisms,  one can use 
the 2-arc diagram in the $XYZ$-coordinate system of $L$ obtained by using 
the following matrix 
\[V=(-\sqrt{2}u_X^t\,-\sqrt{2} u_Y^t\,u_Z^t)=
\left(
\begin{array}{ccc}
0 & 0 & 1\\
1& 1  & 0\\
-1 & 1 & 0
\end{array}
\right)\]
instead of $U$.  By $V$, the 2-arc system $L$ is given by
\[\begin{array}{ll}
a_1&: s(a_1) = (0,0)^0, \quad  t(a_1)=(0,0)^{-1},\\
a_2&: s(a_2)=(1,1)^0, \quad t(a_2)=(-1,1)^0.  
\end{array}\] 
By applying $O(r,\theta)^t$ to this 2-arc system $L$, the $2$-arc diagram 
$D(L)_{21}$ is a 2-arc diagram without crossing, illustrated  in Fig.~\ref{fig:D0arcxy_z}. 
Thus, $\kappa(D(L)_{21})=0$, and  
\[\kappa(L)=\kappa(D(L)_{12})/4=(1/16,1/16,1/16,1/16). \]

\phantom{x}

\phantom{x}

\begin{figure}[hbt]
\centerline{\includegraphics[width=7cm,height=5 cm]{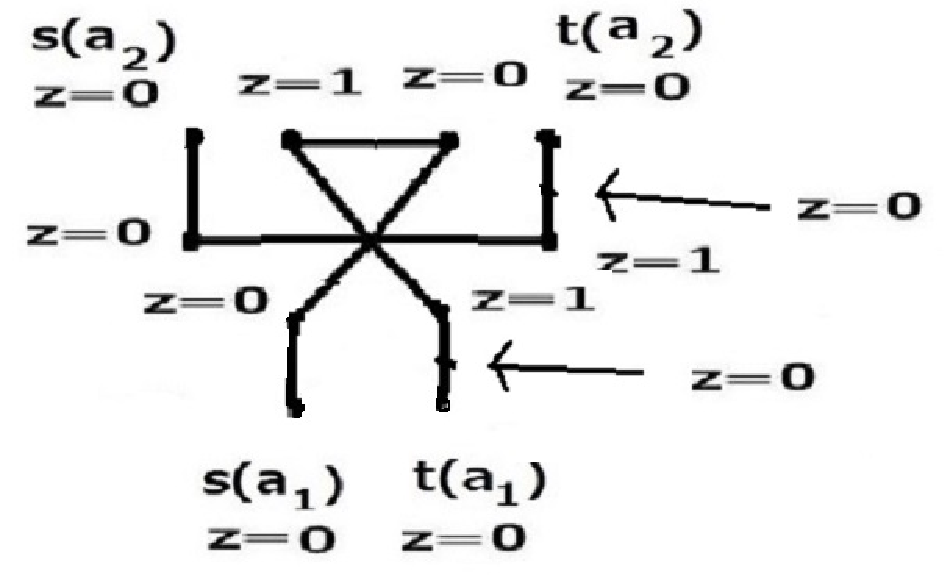}}
\caption{The even $2$-arc system $L=\{a_1, a_2\}$ represented in the 
$xyz$-coordinate system of $R^3$ with the following ordered vertices $(x,y)^z$:\, 
$a_1$: $s(a_1)=(1,0)^0$, $(1,1)^0$, $(2,3)^0$, $(1,3)^1$, $(2,1)^1$, 
$(2,0.5)^0$,  $t(a_1)=(2,0)^0$,\quad  
$a_2$: $s(a_2)=(0,3)^0$, $(0,2)^0$, $(3,2)^1$, $(3, 2.5)^0$, $t(a_2)=(3,3)^0$. }
\label{fig:Larcxy_z}
\end{figure} 

\phantom{x}

\begin{figure}[hbt]
\centerline{\includegraphics[width=13cm,height=9cm]{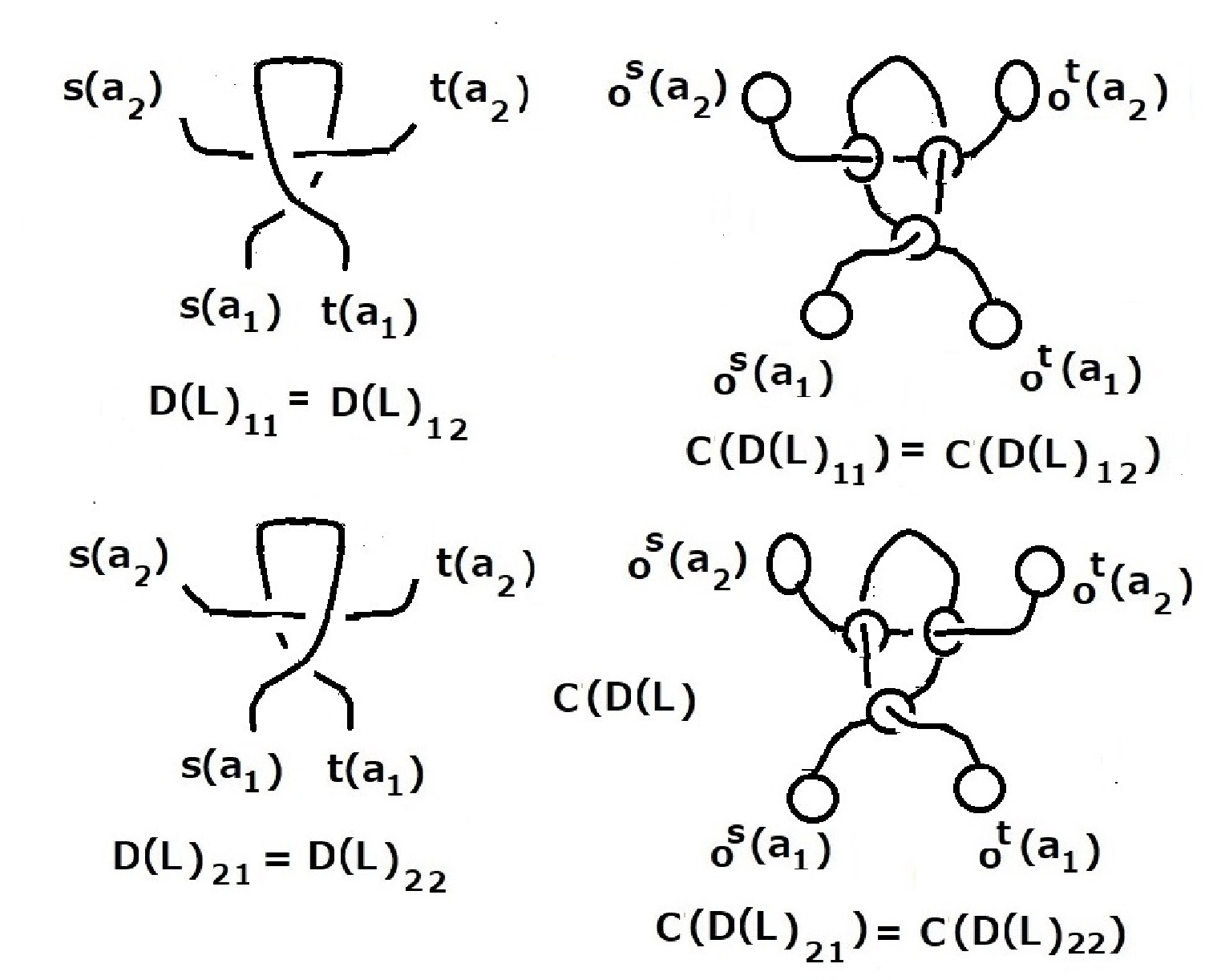}}
\caption{2-arc diagrams $D(L)_{ij}\,(i, j=1,2)$ and the chord diagrams 
$C(D(L)_{ij})\,(i, j=1,2)$}
\label{fig:Darcxy_z}
\end{figure}

\noindent{\bf Example~6.2.}
Let $L$ be an even  $2$-arc system in $R^3$ illustrated  in Fig.~\ref{fig:Larcxy_z}. 
The front pop planes $P_{1j}\,(j=1,2)$ are the plane $z=0$ with ordinal $xy$-orientation 
since $s(a_1)=(1,0,0)$, $u_X=(1\,0\,0)$, $u_Y=(0\,1\,0)$, $u_Z=(0\,0\,1)$. 
The new coordinate system is given by changing the point $(1,0,0)$ into the origin 
$(0,0,0)$, namely by  
the matrix equation $(X\,Y\,Z)=((x\,y\,z)-(1\,0\,0))E$ with the identity matrix $E$. 
The $2$-arc system $L=\{a_1, a_2\}$ is given as follows.
\[\begin{array}{ll}
a_1: s(a_1)=(0,0)^0, (0,1)^0, (1,3)^0, (0,3)^1, (1,1)^1, (1,0.5)^0, t(a_1)=(1,0)^0,\\
a_2: s(a_2)=(-1,3)^0, (-1,2)^0, (2,2)^1, (2,2.5)^0, t(a_2)=(2,3)^0.
\end{array}\]
Let $a_{11}$, $a_{12}$ and $a_{21}$ be oriented edges of $L$ with 
 $a_{12}$, $a_{12}$ in the arc $a_1$ and $a_{21}$ in the arc $a_2$ such that 
\[\begin{array}{lll}
a_{11}&: s(a_{11}) = (0,1)^0,\,  t(a_{11})=(1,3)^0,\\
a_{12}&:s(a_{12}) = (0,3)^1,\, t(a_{12})=(1,1)^1,\\ 
a_{21}&:s(a_{21}) = (-1,2)^0,\, t(a_{21})=(2,2)^1.
\end{array}\]
Let $V(r,0)$ be the matrix  obtained by setting $\theta=0$ to the matrix 
$V(r,\theta)$ used in Example~6.1.
\[V(r,0)=
\left(
\begin{array}{ccc}
1 & 0 & r' \\
0 & 1 & 0 \\
-r' & 0 & 1
\end{array}
\right).\]
In the $XYZ$corodinate system given by $(X\,Y\,Z)=(x\,y\,z) V(r,0)$, 
the oriented edges $a_{11}$, $a_{12}$ and $a_{21}$ are written as follows:
\[\begin{array}{lll}
a_{11}&: s(a_{11}) = (0,1)^0,\,  t(a_{11})=(1,3)^{r'},\\
a_{12}&:s(a_{12}) = (-r',3)^1,\, t(a_{12})=(1-r',1)^{1+r'},\\ 
a_{21}&:s(a_{21}) = (-1,2)^{-r'},\, t(a_{21})=(3-r',2)^{1+3r'}.
\end{array}\]
The 2-arc diagram $D(L)_{1j}\,(j=1,2)$ and the chord diagrams 
$C(D(L)_{1j})\,(j=1,2)$ up to basic equivalences are  illustrated in Fig.~\ref{fig:Darcxy_z}. 
The NT probabilities $\kappa(D(L)_{1j},a_1)$ and  $\kappa(D(L)_{1j},a_2)$ for $j=1,2$  are 
calculated to be
\[\kappa(D(L)_{1j}, a_1)=(1/2, 1/2, 5/8, 7/10),\quad
\kappa(D(L)_{1j}, a_2)=(1/2, 1/2, 3/5, 3/5).\]
Thus,
 \[\kappa(D(L)_{1j})=(1/2, 1/2, 49/80, 13/20)\, (j=1,2).\]

The front pop planes $P_{2j}\,(j=1,2)$ are the plane $z=0$ with opposite $xy$-orientation 
since $s(a_1)=(0,3,0)$, $u_X=(1 \,0 \,0)$, $u_Y=(0\, -1 \,0)$, $y_Z=(0 \,0 \,-1)$.
The new coordinate systems for $P_{2j}\,(j=1,2)$ are given by 
the matrix equation $(X\, Y\, Z)=((x\, y\, z) - (0\, 3\, 0))U$ with the following matrix $U$
\[U=
\left(
\begin{array}{ccc}
1 & 0 & 0 \\
0 & -1 & 0 \\
 0 & 0 & -1
\end{array}
\right).\]
The $2$-arc system $L=\{a_1, a_2\}$ is given as follows.
\[\begin{array}{ll}
a_1: s(a_1)=(2,3)^0, (1,2)^0, (2,0)^0, (1,0)^{-1}, (2,2)^{-1}, (2, 2.5)^0, t(a_1)=(2,3)^0,\\
a_2: s(a_2)=(0,0)^0, (0,1)^0, (3,1)^{-1}, (3, 0.5)^0, t(a_2)=(3,0)^0.
\end{array}\]
Let $a_{11}$, $a_{12}$ and $a_{21}$ be oriented edges of $L$ with $a_{12}$, $a_{12}$ 
in the arc $a_1$ and $a_{21}$ in the arc $a_2$ such that 
\[\begin{array}{lll}
a_{11}&: s(a_{11}) = (1,2)^0,\,  t(a_{11})=(2,0)^0,\\
a_{12}&:s(a_{12}) = (1,0)^{-1},\, t(a_{12})=(2,2)^{-1},\\ 
a_{21}&:s(a_{21}) = (0,1)^0,\, t(a_{21})=(3,1)^{-1}.
\end{array}\]
In the $XYZ$corodinate system given by $(X,Y,Z)=(x,y,z) V(r,0)$, 
the oriented edges $a_{11}$, $a_{12}$ and $a_{21}$ are written as follows:
\[\begin{array}{lll}
a_{11}&: s(a_{11}) = (1,2)^{r'},\,  t(a_{11})=(2,0)^{2r'},\\
a_{12}&:s(a_{12}) =(1+r',0)^{-1+r'},\, t(a_{12})=(2+r',2)^{-1+2r'},\\ 
a_{21}&:s(a_{21}) = (0,1)^0,\, t(a_{21})=(3+r',1)^{-1+3r'}.
\end{array}\]
The 2-arc diagrams $D(L)_{2j}\,(j=1,2)$ are the same inbound 2-arc diagram which is the 
mirror image of the same 2-arc diagram $D(L)_{1j}\,(j=1,2)$, 
as  illustrated in Fig.~\ref{fig:Darcxy_z}. 
By Corollary~5.5,  
the NT probabilities $\kappa(D(L)_{2j})\, (j=1,2)$ are equal to   
$\kappa(D(L)_{1j})\, (j=1,2)$. 
Thus,
 \[\kappa(L) = (1/2, 1/2, 49/80, 13/20). \]

\phantom{x}

\noindent{\bf Acknowledgements.} This work was partly supported by 
JSPS KAKENHI Grant Number JP21H00978 and MEXT Promotion of Distinctive 
Joint Research Center Program JPMXP0723833165 and Osaka Metropolitan 
University Strategic Research Promotion Project (Development of International 
Research Hubs).  

\phantom{x}

\end{document}